\documentclass[preprint, 3p]{elsarticle}

\usepackage{bm}
\usepackage{amsthm}

\usepackage{graphics}
\usepackage{subfigure}
\usepackage{color}
\usepackage{lineno,hyperref}
\usepackage{amsmath}
\usepackage{graphics}
\usepackage{subfigure}
\usepackage{color}
\usepackage{amssymb}
\usepackage{multirow}

\journal{ISTTT (accepted)}

\biboptions{authoryear}

\begin{document}

\begin{frontmatter}

\title{Multi-scale Perimeter Control Approach in a Connected-Vehicle Environment}

\author[1]{Kaidi~Yang}
\ead{kaidi.yang@ivt.baug.ethz.ch}
\author[1,2]{Nan~Zheng\corref{cor1}}
\ead{nan\_zheng@buaa.edu.cn}
\author[1]{Monica~Menendez}
\ead{monica.menendez@ivt.baug.ethz.ch}
\address[1]{Traffic Engineering Group, Institute for Transport Planning and Systems, ETH Zurich, Switzerland}
\address[2]{School of Transportation Science and Engineering, Beihang Univeristy, Beijing, China}
\cortext[cor1]{Corresponding author. Tel.: +41 44 633 32 46}

\begin{abstract}
This paper proposes a novel approach to integrate optimal control of perimeter intersections (i.e. to minimize local delay) into the perimeter control scheme (i.e. to optimize traffic performance at the network level). This is a complex control problem rarely explored in the literature. In particular, modeling the interaction between the network level control and the local level control has not been fully considered. Utilizing the Macroscopic Fundamental Diagram (MFD) as the traffic performance indicator, we formulate a dynamic system model, and design a Model Predictive Control (MPC) based controller coupling two competing control objectives and optimizing the performance at the local and the network level as a whole.  
To solve this highly non-linear optimization problem, we employ an approximation framework, enabling the optimal solution of this large-scale problem to be feasible and efficient. Numerical analysis shows that by applying the proposed controller, the protected network can operate around the desired state as expressed by the MFD, while the total delay at the perimeter is minimized as well. Moreover, the paper sheds light on the robustness of the proposed controller. This multi-scale hybrid controller is further extended to a stochastic MPC scheme, where connected vehicles (CV) serve as the only data source. Hence, low penetration rates of CVs lead to strong noises in the controller. This is a first attempt to develop a network-level traffic control methodology by using the emerging CV technology. We consider the stochasticity in traffic state estimation and  the shape of the MFD. Simulation analysis demonstrates the robustness of the proposed stochastic controller, showing that efficient controllers can indeed be designed with this newly-spread vehicle technology even in the absence of other data collection schemes (e.g. loop detectors).
\end{abstract}

\end{frontmatter}

\section{Introduction}
\label{sec:1}

Real-time traffic control strategies at the aggregated level have been receiving significant research attention, in particular the perimeter control (also known as gating). The basic concept is to restrict the incoming flow through traffic signals at the boundaries of a pre-defined region to prevent congestion inside. The controllers of this type are typically designed based on the Macroscopic Fundamental Diagram (MFD), also known in the literature as the Network Fundamental Diagram (NFD). Examples of such control strategies can be found in  \citet{Keyvan-Ekbatani2012,geroliminis2013optimal,
haddad2013cooperative,aboudolas2013perimeter,
haddad2015robust,ramezani2015dynamics}. 

The concept of an MFD or NFD was initially proposed in \citet{godfrey1969} and followed in \citet{mahmassani1987performance} and \citet{daganzo2007}. The demonstration of the existence of the MFD with dynamic features from field data was firstly reported in \citet{geroliminis2008existence}, showing that traffic in urban single-mode regions exhibits an aggregated relationship between the space-mean network flow (traffic throughput) and the network density (traffic accumulation). During the past decade, plenty of research efforts have been devoted to further investigate the properties and develop approximations of the MFD model. An interested reader could refer to \citet{Yildirimoglu2014} and \citet{Leclercq2014Macroscopic} for a review of recent developments in MFD. For the control research field, the most important contribution of the MFD model lies in the fact that the model features the so-called critical accumulation of the network, which serves as the control goal and facilitates the design of various control schemes for traffic signals in real time. It has been shown in the literature that MFD-based control strategies are effective in regulating the global traffic performance; examples for single-mode networks can be found in \citet{haddad2012stability}, \citet{Keyvan-Ekbatani2012,Keyvan-Ekbatani2015}, and \citet{haddad2015robust}.  

As the network-wide control has increasingly gained relevance, alternative and more complex concepts have been developed. For example, approaches have been proposed to partition heterogeneous networks (where congestion levels are unevenly distributed) into a small number of homogeneous regions, and to apply perimeter control to the inter-regional flows along the boundaries between those regions. This is known as multi-region control \citep[see][]{aboudolas2013perimeter,
ramezani2015dynamics,kouvelas2016enhancing}. The inter-transferring flows are controlled at the intersections located at the border between regions, so as to distribute traffic in an optimal way and minimize the total travel costs of all regions. This type of control can be viewed as a high-level regional control scheme and may be combined with other strategies (e.g. local or distributed controllers) in a hierarchical control framework \citep{kouvelas2016enhancing}. Recent advances towards this direction are reported by \cite{daganzo2007}, \cite{Keyvan-Ekbatani2012}, \cite{geroliminis2013optimal}, \citet{aboudolas2013perimeter}, \cite{ramezani2015dynamics}, and \cite{kouvelas2016enhancing}. Furthermore, some studies have extended the network-wide control to multimodal networks. \citet{ampountolas2014perimeter} developed perimeter control strategies for a bimodal network, recognizing buses have larger transport efficiency than cars,  therefore certain amount of buses should operate in the system to optimize passenger throughput. 
\citet{chiabaut2015evaluation} looked at the MFD model of arterial roads from a passenger flow perspective, and developed optimization strategies for bus operations. These last two references provide approaches for optimizing the passenger mobility through real-time control.

Regardless of what type of algorithms are applied, the optimal flow allowed to enter the network has to be distributed to the local intersections. In other words, as long as the controllers determine the total flow, the flow should be allocated to each individual perimeter intersection through the control of the traffic signals. Nearly all the existing studies have emphasized that delay may be caused at local intersections when applying perimeter control \citep{Keyvan-Ekbatani2012,geroliminis2013optimal,haddad2013cooperative,aboudolas2013perimeter,
hajiahmadi2015optimal,haddad2017optimal}. However, to the best of our knowledge, no work has quantitatively and systematically treated the delay at these intersections when designing the controllers for network-wide applications. Considering detailed performance of the local intersections may complicate the dynamics of the system, and consequently the optimization problem. Thus, given its challenging complexity, controllers treating specifically the intersections are rarely reported in the literature. Although some initial efforts have been made (e.g. combining adaptive traffic signal settings inside the network \citep{kouvelas2016enhancing}, 
and considering the queue length at the perimeter \citep{keyvan2016combination,haddad2017optimal}, the interaction between the network level perimeter control and the local level intersection control has not been fully considered.

Motivated by the discussion above, this paper proposes a novel perimeter control approach that aims at improving the global traffic performance while simultaneously accounting for new technologies in transportation. 
Specifically, we propose and design control algorithms that work efficiently in a connected vehicle environment.  
The connected vehicle technology has been attracting increasing attention in the traffic control field, thanks to its capability of detailed and anticipative information provision. 
This information can be used not only to measure the current traffic situation, but also to predict traffic states in the future. 
This facilitates the application of the Model Predictive Control \citep{garcia1989model} (MPC), which takes the predicted information and current states as inputs to optimize the current control actions while taking into consideration the future performance. Given its ability to handle complex and dynamic system models, MPC-based approaches have been studied for traffic control \citep[see for example,][]{geroliminis2013optimal,hajiahmadi2015optimal}. 

The contributions of this paper are two-fold. (1) We develop a multi-scale control algorithm that optimizes traffic performance at the network level, and at the local level (i.e. perimeter intersections). This is challenging, as two competing control objectives are coupled into an integrated multi-scale control scheme. We will show that by applying the proposed multi-scale controller, the protected network can operate around the desired state as expressed by the MFD, while the total delay at the perimeter is minimized as well. (2) We apply the proposed controller in a connected-vehicle environment, for which the robustness of the control is enhanced.  Although connected vehicle technology has attracted much attention in managing simple intersections \citep{guler2014using,yang2016isolated}, this work, to the best of our knowledge, is the first attempt to develop a network-level methodology using such technology. 
Evidently, using also traditional data sources such as loop detectors, video cameras and probe vehicles would be helpful to measure traffic accumulations \citep[e.g. see][]{ambuhl2016data} for control purposes. However, loop detectors are only expensive but also inflexible in location which results in biased estimation, while video cameras are be sensitive to weather and converage conditions. In comparison, the data provided by connected vehicles is richer, more diverse and detailed. These advantages make CV environment exclusively suitable for developing multi-scale control. On the other hand, the application of connected vehicle technology is expected to keep growing rapidly in most places around the world. The data will have a much wider coverage of the road network. The challenge, nevertheless, lies in the uncertainties due to its limited penetration rates, especially in the near future. In perimeter flow control literature, control schemes are developed taking into account uncertainties in the system.   
Some recent works proposed robust perimeter control strategies for unimodal networks \citep{geroliminis2013optimal,haddad2015robust}, and multimodal networks \citep{ampountolas2014perimeter}. These works assumed to have perfect information on the accumulation of vehicles, thus may not be applicable to the connected vehicle environment. We consider connected vehicles as the only source for obtaining traffic information, which brings noise into the measurement of vehicle accumulations. Previous researches have already shown that low penetration rates of probe vehicles reduce the accuracy of vehicle density estimation \citep[see][]{gayah2013using,nagle2014accuracy,nagle2015comparing,
du2016deriving, ambuhl2016data}.
 In a fully connected vehicle environment, a biased state estimation would have significant impact on control. Hence, we will present an effective stochastic control approach that can improve the design of control strategies utilizing information solely coming from this newly-spread vehicle technology.

The rest of the paper is organized as follows. Section \ref{sec:2} introduces the general framework of a dynamic traffic system model and the control problem. Section 3 describes in detail the design of the proposed control algorithms: a multi-scale perimeter flow controller and a stochastic controller that takes into account system uncertainties brought by the limited information from connected vehicles. Section 4 analyzes the performance of the proposed control scheme and carries out a comparative analysis to the classical control schemes. Section 5 further investigates the robustness of the proposed control scheme and demonstrates the importance of the multi-scale treatment. Concluding remarks and discussions are given in Section 6.

\section{General framework}
\label{sec:2}
In this section, we present the general framework for the dynamic traffic system model, and the control problem under consideration. For the reader's convenience, Table \ref{tab:nomen} provides the notation for the main variables and parameters in this paper.

Consider a typical urban city of single-center structure, with two regions. Let us denote the center region as ``1'' and the outside region as ``2'', where the perimeter between the two regions is the perimeter where we apply the control via traffic signals. The two levels of control under consideration are (i) the city center network which attracts large demand (i.e. network level); and (ii) the intersections at the perimeter, where the perimeter control is implemented (i.e. local level).
At the network level, the dynamics of the system are described by an MFD $O=G(n)$ where the aggregated traffic completion flow $O$ (outflow, vehs per time unit) is a function $G(\cdot)$ of the traffic accumulation $n$ (number of vehicles, vehs). The MFD is assumed to be known for the given network.

The change in traffic accumulations of the region, reflecting the network-level traffic state, can be represented by the evolution of the accumulation $n$, which is captured by mass conservation equations without the need for detailed traffic information (such as routing at link level) within the network. 
Eq.(\ref{eq:n11c})-Eq.(\ref{eq:n12c}) show the time-discretized dynamics at the network level. The cycle length is denoted as $C$ [hr]. 
\begin{align}
n_{11}(k+1)&= n_{11}(k) + D_{11}(k)C + \beta_{21}(k)C - \frac{n_{11}(k)}{n_{11}(k)+n_{12}(k)}G(n_{11}(k)+n_{12}(k))C \label{eq:n11c}\\
n_{12}(k+1)&= n_{12}(k)+ D_{12}(k)C - \beta_{12}(k)C \label{eq:n12c}
\end{align}

\begin{table}[htbp]
\centering
\caption{Nomenclature}
\centering
\begin{tabular}{c|p{13cm}}
\hline \hline
region $1$ & center region of the network\\\hline
region $2$ & periphery region of the network  \\\hline
$I$ & set of intersections \\\hline
$M^i$ & set of streams of an intersection $i$, $i\in I$ \\\hline
$P^i$ & set of phases of intersection $i$, $i \in I$ \\\hline
$C$ & cycle length of each intersection \\\hline
$k$ & index of the current cycle \\\hline
$L$ & prediction horizon (number of cycles) \\\hline
$l$ & index of predicted cycle from the current cycle\\\hline
$G(\cdot)$ & MFD that relates the completion flow to the traffic accumulation in a network \\\hline
$v$ & absolute value of the slope of the left-branch of the MFD (free flow states) \\\hline
$w$ & absolute value of the slope of the right-branch of the MFD (congested states) \\\hline
$n_{cr}$ & critical accumulation in the MFD \\\hline
$\hat{n}_{ab}(k)$ & measured traffic accumulation for vehicles in Region $a$ with destination in Region $b$ at cycle $k$, $a,b=1,2$. It is the initial traffic accumulation for the proposed MPC model \\\hline
$\hat{x}_m^i(k)$ & measured queue length of stream $m$ at intersection $i$ during cycle $k$, $i\in I$ , $m \in M^i$. $\hat{x}_m^i(k)$ is considered as the initial value for the proposed MPC model \\\hline
$n_{ab}(k+l|k)$ & predicted accumulation of vehicles in Region $a$ with  destination in Region $b$ at cycle $k+l$  based on information available at cycle $k$, $0\leq l\leq L$, $a,b=1,2$ \\\hline
$D_{ab}(k+l|k)$ &  predicted newly generated demand, generated from region $a$ with destination of region $b$ at cycle $k+l$  based on information available at cycle $k$, $0\leq l\leq L$, $a,b=1,2$ \\ \hline
$\beta_{ab}(k+l|k)$ & predicted controlled flow transferring from Region $a$ to Region $b$ at cycle $k+l$ based on information available at cycle $k$, $0\leq l\leq L$, $a,b=1,2$\\\hline
$s_{mp}^i$ & maximum discharging flow of stream $m$ in cycle $p$ at intersection $i$, $i \in I$, $m\in M^i$, $p\in P^i$ \\\hline
$\mu_m^i(k+l|k)$ & predicted departure flow of stream $m$ at intersection $i$ in cycle $k+l$ based on information available at cycle $k$, $0\leq l\leq L$, $i \in I$, $m\in M^i$ \\\hline
$q_m^i(k+l|k)$ & predicted newly generated arrival flow of stream $m$ at intersection $i$ in cycle $k+l$ based on information available at cycle $k$, $0\leq l\leq L$, $i \in I$, $m\in M^i$ \\\hline
$x_m^i(k+l|k)$ & predicted queue length of stream $m$ at intersection $i$ in cycle $k+l$ based on information available at cycle $k$, $0\leq l\leq L$, $i \in I$, $m\in M^i$.  \\\hline
$g_p^i(k+l)$ & green time ratio of phase $p$ at intersection $i$ at cycle $k+l$, $0\leq l\leq L$, $i \in I$, $p\in P^i$ \\\hline
$\alpha_m^i$ & percentage of outflow vehicles on stream $m$ at intersection $i$, $i \in I$, $m\in M^i$ \\\hline
$g_{\max}^i$ & maximum allowed total green time ratio at intersection $i$,  $i \in I$ \\\hline
$g_{p,\min}^i$ & minimum allowed green time ratio for phase $p$ at intersection $i$,  $i \in I,~ p \in P^i$ \\\hline
\hline
\end{tabular}
\label{tab:nomen}
\end{table}

In the equations, $n_{ab}(k)$ [veh] represents the traffic accumulation in region $a$ with destination in region $b$ at cycle $k$ (independently of the origin of the individual trips). Hence, $n_{11}(k)$ denotes the traffic accumulation in the city center with a destination in the city center; $n_{12}(k)$ is the traffic accumulation in the city center which will end trips in the outside region; likewise, $n_{21}(k)$ is the traffic accumulation in the outside region that will end their trips in the city center.  $D_{ab}(k)$ [veh/hr] denotes the newly generated demand from region $a$, with a destination to region $b$ at cycle $k$. Therefore $D_{11}(k)$ is the internal demand having destinations inside the city center; $D_{21}(k)$ is the demand coming from outside while ending their trips inside of the center center; $D_{12}(k)$ is the demand generated in the center region with destination to the outside region; and $D_{22}(k)$ is not considered here. $\beta_{21}(k)$ and $\beta_{12}(k)$ [veh/hr] are the controlled flow, entering and leaving the center region, respectively (these are important variables to the proposed control scheme). The last term in Eq.(\ref{eq:n11c}) represents the uncontrolled trip completion flow inside the city center. Since the MFD function $G$ is known and the demand $D_{ab}(k)$ can be learned  without requiring detailed Origin-Destination (OD) information, traffic accumulation Eq.(\ref{eq:n11c})-Eq.(\ref{eq:n12c}) can be monitored in real time. Uncertainties in the demand $D_{ab}(k)$ will be discussed in Section \ref{sec:5}. 

Initially, we assume that the network is equipped with a sufficient number of connected vehicles (i.e. high penetration rate) so that we can properly measure the states (i.e. traffic accumulation and queue lengths) of the network. This assumption is relaxed in Section \ref{sec:5}, and its impact on the traffic performance will be discussed later. Theoretically, the method proposed in this paper can also use the data from other sources such as video cameras, probe vehicles, and historical data. However, in practice, the method is expected to work with better accuracy with connected vehicle data, as we assume real-time information is available on the dynamic origin-destination (OD) at the network level and the arrival flow at the intersection level (a priori knowledge on which intersection each vehicle uses to enter the city center). 

For the local level, we track the dynamics of the queues at the perimeter intersections. Denote the set of the intersections as $I$, the set of streams for each intersection $i$ as $M^i$, and the phases for each intersection as $P^i$. We assume a constant cycle length for all intersections $C$ [hr]. This assumption reduces complexity in formulating the control problem, allowing an extensive numerical analysis with moderate computational burden. Alternatively, cycle length and other signal settings can be dynamically formulated, and easily integrated into the proposed model.    

We refer to the departure flow for each direction $m$ of intersection $i$ at cycle $k$ as $\mu_m^i(k)$ [veh/hr]. 
$\mu_m^i(k)$ is estimated based on the arrival flow $q_m^i(k)$ [veh/hr], the queue length $x_m^i(k)$ [veh], the discharging flow $s_{mp}^i$ [veh/hr], the cycle length $C$ [hr], and the allocated green time ratio $g_p^i(k)$ [-] \citep{liu2008reverse,zhao2011fast,li2013optimality}.
For each intersection $i$, the maximum discharging flow rate $s_{mp}^i$ is a time-invariant variable which takes one of the following two values: 1) the saturation flow rate, if green signal is given to stream $m$ in phase $p$ at intersection $i$; 2) 0, otherwise.  $s_{mp}^i$ only represents the ability of stream $m$ at intersection $i$ to discharge vehicles in phase $p$. It is independent of the arrival flow of vehicles. The actual departure rate $\mu_m^i(k)$ is, on the other hand, bounded by the demand at this stream.
Note that we assume that the arrival flow $q_m^i(k)$ [veh/hr], the queue length $x_m^i(k)$ [veh], the maximum discharging flow $s_{mp}^i$ [veh/hr], and the allocated green time ratio $g_p^i(k)$ can have different values among intersections. 
The evolution of the queue length $x_m^i(k+1)$ can be captured by the mass conservation law, as well, i.e. Eq.(\ref{eq:queueCE}). 
\begin{align}
x_m^i(k+1) = x_m^i(k) + q_m^i(k)C - \mu_m^i(k)C \label{eq:queueCE}
\end{align}

For the proposed framework at the local level, we take into account not only the directions for accommodating flows  $\beta_{21}(k)$ and $\beta_{12}(k)$, but also the general directions of a typical signalized intersection. 
This is important to bridge the network level flow physically to the local intersections. The existing works \citep{haddad2012stability, geroliminis2013optimal} rarely consider detailed traffic signal settings, focusing solely on the network inflow and outflow that influence the evolution of the total accumulation. 
The rest of the traffic streams, e.g. flows parallel to the boundaries, are neglected. Such treatment may be sufficient for proving the concept of a network-level control approach. As shown by the dynamic equations and the control problem to solve in the later text, ignoring the other directions over-simplifies the formulation of the control problem and overlooks the local queues which possibly negate some of the benefits of the perimeter control. 
Other studies, such as \citet{aboudolas2013perimeter, Keyvan-Ekbatani2015, kouvelas2016enhancing} optimize traffic signals at perimeter or boundaries between regions, utilizing micro-simulations where it is possible to change the signal control plan for each individual intersection. Nevertheless, no detailed methodology has been reported nor formulated on the optimization of signal timing allocations at the local level. The proposed framework in this paper aims to fill this gap.

\section{Methodology for the multi-scale control}
\subsection{An MPC approach for multi-scale control}
\label{sec:3.1}
An MPC approach is proposed to the multi-scale controller which integrates the network level perimeter control and the local level intersection control through detailed signal optimization. 
The proposed MPC relies on a dynamic traffic model that couples both the network and the intersections. 
It calculates the green time ratios at each intersection (i.e. control variables) to optimize the total travel cost in the current signal cycle, while taking into consideration the costs in the future cycles. 
This is achieved by solving an optimization problem with a finite-time horizon ($N$ cycles), where we obtain the evolution of traffic accumulations and the predicted optimal green ratios for each cycle, but only execute the green ratios for the current cycle. These features make the proposed MPC approach fundamentally different, while more complicated to solve than the ones in previous studies.

There are two types of inputs to the MPC: the current traffic states (i.e. traffic accumulation and queue lengths at each intersection), and the traffic demands.
The traffic accumulation and queue lengths are measured using information communicated from the connected vehicles. The traffic demands and the region-level ODs can be predicted in real time through the CVs alone or combined with historical data.

 \begin{figure}[t]
\centering
\includegraphics[width=0.9\textwidth]{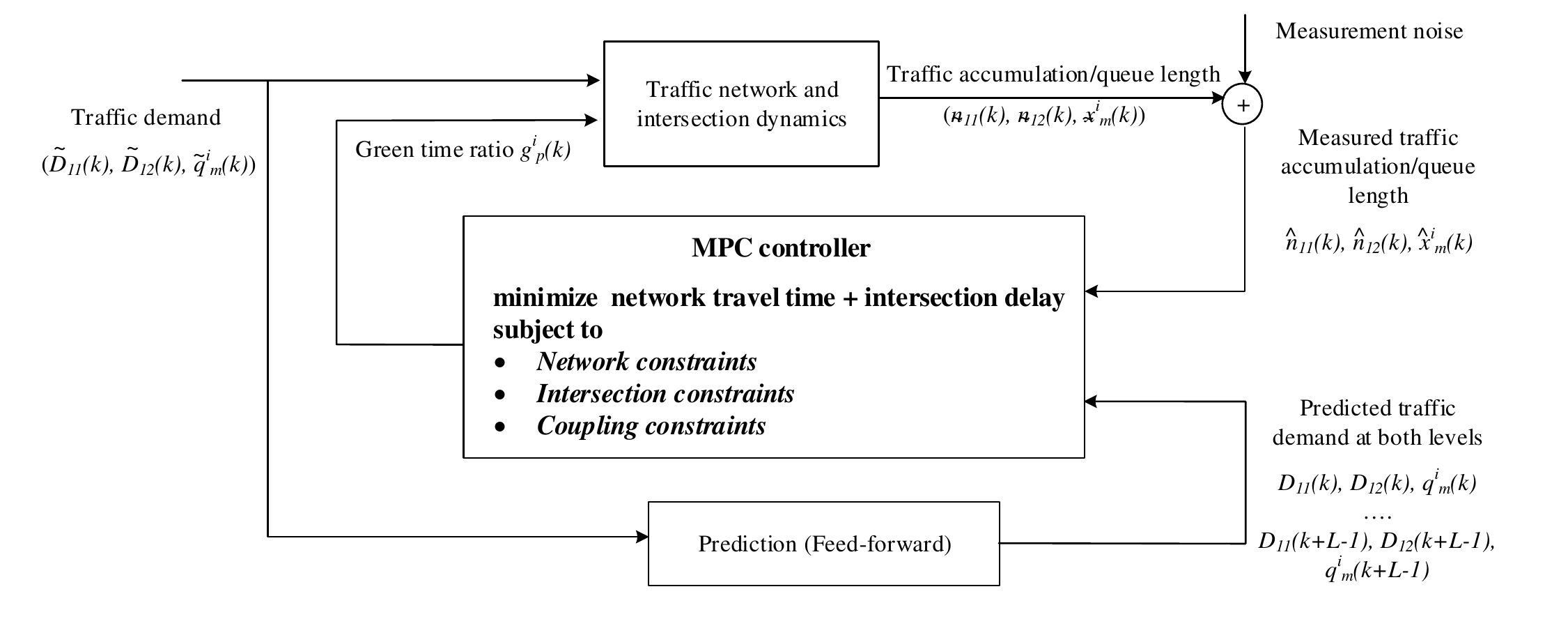}
\centering
\caption{Control diagram for MPC.}
\label{fig:robust}
\end{figure}

The control diagram in Fig.\ref{fig:robust} illustrates graphically the control process. Note that we separate the real system dynamics (at both network and local level) and the MPC prediction models in Fig.\ref{fig:robust}. The MPC prediction model can be different from the system dynamics. 
$g_p^i(k)$ are the control variables (i.e., output of the MPC controller). 
$\tilde{D}_{ab}(k)$ and $\tilde{q}_m^i(k)$ represent the current real demand at the network and local level at cycle $k$, respectively; 
$D_{ab}(k+l)$ and $q_m^i(k+l)$, $0\leq l \leq L-1$ represent the predicted demand at the network level and the intersection level, respectively.   
$\tilde{n}_{ab}(k)$ and $\tilde{x}_{m}^i(k)$ represent the current real traffic accumulation and queue length in the traffic system, respectively;  
$\hat{n}_{ab}(k)$ and $\hat{x}_{m}^i(k)$ represent the measured traffic accumulation and queue length at cycle $k$, respectively, both are considered to be obtained through connected vehicles. 
Notice that no measurement noise nor prediction error is incorporated at this stage, i.e. $D_{ab}(k)=\tilde{D}_{ab}(k)$, $q_m^i(k) = \tilde{q}_m^i(k)$, $\hat{n}_{ab}(k) = \tilde{n}_{ab}(k)$, $\hat{x}_m^i(k) = \tilde{x}_m^i(k)$. This assumption will be relaxed later where additional treatment is integrated. 

At each cycle $k$, we solve the MPC optimization model over a moving time horizon of the next $L$ cycles. 
The detailed formulation of the optimization problem is displayed in Eq.(\ref{eq:obj})-(\ref{eq:xinit}).
The inputs to the model are (a) the predicted demand in the next $L$ cycles at both the network and local level, (b) the measured traffic accumulation, and (c) the queue length. 
The output of the model are the optimal green time ratios for each phase at each intersection during the current cycle, which serve as a feedback to the traffic system.

The objective function Eq.(\ref{eq:obj}) minimizes the total travel cost $J_D$ of the whole system, incurred at the perimeter intersections and within the controlled network subject to constraints Eq.(\ref{eq:n11d})-Eq.(\ref{eq:xinit}). 
The first two terms in Eq.(\ref{eq:obj}), multiplied by the cycle length, represent the network travel time in the city network during a given cycle; the third term (also multiplied by the cycle length) represents the total travel delay at the perimeter intersections also during a given cycle. 
This is computed for all future $L$ cycles (the horizon). Notice that the last term excludes the outflow streams from the network $M^i_{\rm{out}}$ (explanation on the exclusion is given below). 
In Eq.(\ref{eq:obj}), the two objective criteria (travel time in the city network and delay at the perimeter intersections) have the same weights. However, the weights can be adjusted for different policy-oriented purposes.

\begin{align}
\min~&J_D = C\sum_{l=1}^{L}{\Big(n_{11}(k+l|k)+n_{12}(k+l|k)+\sum_{i\in I}{\sum_{m\in M\backslash M_{\rm{out}}^i}{x_{m}^i(k+l|k)}}\Big)} \label{eq:obj}\\
\rm{s.t.}~&n_{11}(k+l+1|k) = n_{11}(k+l|k) + D_{11}(k+l|k)C + \beta_{21}(k+l|k)C  \notag \\
&~~~~~~~-\frac{n_{11}(k+l|k)}{n_{11}(k+l|k)+n_{12}(k+l|k)}G\Big(n_{11}(k+l|k)+n_{12}(k+l|k)\Big)C,~\forall 0\leq l\leq L-1 \label{eq:n11d}\\
&n_{12}(k+l+1|k) = n_{12}(k+l|k) + D_{12}(k+l|k)C - \beta_{12}(k+l|k)C,~\forall 0\leq l\leq L-1 \label{eq:n12d}\\
&\beta_{21}(k+l|k) = \sum_{i \in I}{\sum_{m \in M_{\rm{in}}^i}{\mu_m^i(k+l|k)}},~\forall 0\leq l\leq L-1 \label{eq:betamu21}\\
&\beta_{12}(k+l|k) = \sum_{i \in I}{\sum_{m \in M_{\rm{out}}^i}{\mu_m^i(k+l|k)}},~\forall 0\leq l\leq L-1 \label{eq:betamu12}\\
& 
\mu_m^i(k+l|k) =  \min\Big\{\frac{\alpha_m^i}{|I|}\frac{n_{12}(k+l|k)}{n_{11}(k+l|k)+n_{12}(k+l|k)}G\Big(n_{11}(k+l|k)+n_{12}(k+l|k)\Big), \sum_{p\in P^i}{s_{mp}^ig_p^i(k+l)}\Big\},\notag \\
&~~~~~~~~~~~~~~~~~~~~~~~  \forall m \in M^i_{\rm{out}},~\forall i \in I, ~\forall 0\leq l\leq L \label{eq:muout}\\
&\mu_m^i(k+l|k) = \min\{x^i_m(k+l|k)/C  + q^i_m(k+l|k), \sum_{p\in P^i}{s_{mp}^ig_p^i(k+l)}\},~\forall m \in M^i\backslash M^i_{\rm{out}},~\forall i \in I,~\forall 0\leq l\leq L-1\label{eq:muother}\\
&x^i_m(k+l+1|k) = x^i_m(k+l|k) + q^i_m(k+l|k)C - \mu_m^i(k+l|k)C,~\forall m \in M^i\backslash M^i_{\rm{out}},~\forall i \in I,~\forall 0\leq l\leq L-1 \label{eq:xother}\\
&\sum_{p\in P}{g_p^i(k+l)} \leq g_{\max}^i,~\forall i \in I, ~0\leq l\leq L-1 \label{eq:gallow}\\
&g_p^i(k+l) \geq g_{p,\min}^i,~\forall p \in P^i, \forall i \in I,\forall 0\leq l\leq L-1\label{eq:gmin}\\
& n_{1j}(k|k) = \hat{n}_{1j}(k),~j=1,2
\label{eq:n11init}\\
& x_m^i(k|k) = \hat{x}_{m}^i(k), ~\forall m \in M^i\backslash M^i_{\rm{out}},~\forall i \in I \label{eq:xinit}
\end{align}

Constraints Eq.(\ref{eq:n11d}) and Eq.(\ref{eq:n12d}) represent the mass conservation of vehicle accumulations at the network level (they are modelled after Eq.(\ref{eq:n11c}) and Eq.(\ref{eq:n12c})). 
Eq.(\ref{eq:betamu21}) and Eq.(\ref{eq:betamu12}) are ``coupling'' constraints which link the network and the local level. 
Here $M_{\rm{in}}^i$ and $M_{\rm{out}}^i$ are two disjoint subsets of set $M$, which represent the set of inflow streams to the network and the set of outflow streams from the network at intersection $i$, respectively. 
Eq.(\ref{eq:betamu21}) reflects that  the vehicles departing from streams in $M_{\rm{in}}^i$ are the inflow entering into region 1. Likewise, Eq. (\ref{eq:betamu12}) holds because vehicles departing from streams in  $M_{\rm{out}}^i$ are the outflow leaving from region 1. 
Eq.(\ref{eq:muout}) and Eq.(\ref{eq:muother}) ensure that the departure flow of each stream is taken as the minimum of the demand and the capacity determined by the allocated green time.
The first term in Eq.(\ref{eq:muout}) represents the ``demand'' leaving the network. It is the flow arriving at the outflow streams $M_{\rm{out}}^i$ of each intersection, the sum of which is determined by the outflow of the network.
In the prediction model, we assume that the outflow of the network is evenly distributed across all intersections. 
For each intersection, there might exist multiple outflow streams accommodating the flows that leave the network. 
$\alpha_m^i$ represents the percentage of outflow vehicles on each of them. The value of $\alpha_m^i$ can be dynamically updated by online estimation (e.g. a Kalman filter) from CV data while integrating historical information. Note that obtaining $\alpha_m^i$ from traditional data sources (such as fixed loop detectors) may require more effort while resulting in less accuracy.
The first term in Eq.(\ref{eq:muother}) is the total demand for each stream at each intersection, determined jointly by the queue length and the arrival flows. 
Here, $q_m^i(k+l|k)$ represents the newly arriving flow for stream $m$ at cycle $k+l$ at intersection $i$, predicted at cycle $k$. 
The second terms in both Eq.(\ref{eq:muout}) and Eq.(\ref{eq:muother}) represent the capacity, determined by the allocated green time ratios. 
Eq.(\ref{eq:xother}) describes the queue  dynamics at the local level using basic conservation law. 
It does not include the outflow streams of the network, as 1) the total arrival to these streams is decided by the outflow of the network; 2) in the objective function Eq.(\ref{eq:obj}), these outflow vehicles are already included in  $n_{12}(t)$. 
 Eq.(\ref{eq:gallow}) and Eq.(\ref{eq:gmin}) are the physical constraints imposed on the green time ratio. $g_{\max}^i$ is the maximum allowed total green ratio across all phases, defined as $1-\eta/C$ where $\eta$ is the total lost time. $g_{p,\min}^i$ is the minimum duration of the green signal sufficient to discharge a given number of vehicles (e.g. 2 vehicles). 
Both $g_{\max}^i$ and $g_{p,\min}^i$ are given variables, which can be determined by the configurations of the intersections. 
Constraints Eq.(\ref{eq:n11init}) and Eq.(\ref{eq:xinit}) define the initial traffic accumulations, $n_{1j,init}$, and the initial queue lengths, $x_{m,init}^i$. 

As $n_{11}(k+l|k)$, $n_{12}(k+l|k)$, $\beta_{12}(k+l|k)$, $\beta_{21}(k+l|k)$, $\mu_m^i(k+l|k)$, $x_m^i(k+l|k)$ are all functions of $g_p^i(k+l)$, the only decision variables are the green time ratios $g_p^i(k+l)$ for $p \in P^i,~i \in I,~\forall 0\leq l \leq L$. 

It can been seen that the optimization problem defined by  Eq.(\ref{eq:obj})-(\ref{eq:xinit}) is highly complex. To reduce the non-convexity and facilitate the solution searching, constraints Eq.(\ref{eq:n11d}) and Eq.(\ref{eq:muout}) are converted into a piecewise linear function. Details are presented in Section \ref{sec:3.2}.

\subsection{Linearization of the model}
\label{sec:3.2}
In this subsection, we reformulate the optimization problem Eq.(\ref{eq:obj})-Eq.(\ref{eq:xinit}) through linearization techniques. For presentation simplicity, we neglect the cycle index $k$ in this subsection.

A triangular-shape MFD  \citep[used also in ][]{haddad2012stability, haddad2013cooperative}） is employed in the prediction model. Note that for the numerical simulation, an empirical MFD (with scatters and noises) is utilized. This makes the prediction model different from the simulation model (used later in Section 4, Section 5.2, and Section 5.3). This is an important treatment for the MPC-based approach. Also note that choosing alternative MFD shapes  for the prediction model  is straightforward in our framework. Our results show, however, that the triangular MFD represents a good compromise between accuracy and efficiency, and works sufficiently well. 

The MFD is represented as 
\begin{align}
G(n) = \min\{vn,(w+v)n_{cr}-wn\}
\end{align}
where $v,~w>0$ are the absolute value of the slopes of the left and the right branch of the MFD, respectively; $n_{cr}$ is the critical vehicle accumulation.  

The proposed MPC with its current form Eq.(\ref{eq:obj})-Eq.(\ref{eq:xinit}) is highly nonlinear and nonconvex due to the following terms in Eq.(\ref{eq:n11d}) and Eq.(\ref{eq:muout}), respectively. 
\begin{align}
O_{1j} &= \frac{n_{1j}}{n_{11}+n_{12}}G(n_{11}+n_{12}), j = 1,2
\end{align}
Hence, we rewrite $O_{1j}, j = 1,2$ as follows. 
\begin{align}
O_{1j} &= \min\Big\{vn_{1j},\frac{(w+v)n_{cr}n_{1j}}{n_{11}+n_{12}}-wn_{1j}\Big\},~j = 1,2 \label{eq:o1}
\end{align}

It can be seen from Eq.(\ref{eq:o1}) that the only nonlinear terms are $\frac{n_{1j}}{n_{11}+n_{12}}, j = 1,2$. 
To solve the problem efficiently, we apply linearization technique to this term with a first-order Taylor polynomial approximation. 
Denote the initial accumulation as $(\hat{n}_{11}, \hat{n}_{12})$. Specifically, we linearize the second term of $O_{11}$ at $(\hat{n}_{11}, n_{cr}-\hat{n}_{11})$ and the second term of $O_{12}$ at $(n_{cr}-\hat{n}_{12}, \hat{n}_{12})$. With simple calculations, we can then rewrite $O_{11}$ and $O_{12}$ as follows. 
\begin{align}
O_{1j} \approx \bar{O}_{1j} = \min\Big\{vn_{1j},\gamma_{j0} +\gamma_{j1}n_{11}+\gamma_{j2}n_{12}\},~j = 1,2\label{eq:lino11}
\end{align}
where 
\begin{align}
\gamma_{10} &= (w+v)\hat{n}_{11},~
\gamma_{11} = v-(w+v)\hat{n}_{11}/n_{cr},~
\gamma_{12} = -(w+v)\hat{n}_{11}/n_{cr} \label{eq:gamma1}\\
\gamma_{20} &= (w+v)\hat{n}_{12},~
\gamma_{21} = -(w+v)\hat{n}_{12}/n_{cr},~
\gamma_{22} = v-(w+v)\hat{n}_{12}/n_{cr},~
\label{eq:gamma2}
\end{align}

Notice that we linearize $O_{11}$ and $O_{12}$ at $(\hat{n}_{11}, n_{cr}-\hat{n}_{11})$ and $(n_{cr}-\hat{n}_{12}, \hat{n}_{12})$ instead of the original initial accumulation $(\hat{n}_{11}, \hat{n}_{12})$. The reason is summarized in Lemma 1.
 
\textbf{Lemma 1.} By linearizing $O_{1j}$  as $\bar{O}_{1j}$ (Eq.(\ref{eq:lino11}), we obtain 
\begin{align}
\bar{O}_{1j} = O_{1j}, \mbox{ if } n_{11}(k)+n_{12}(k) \leq n_{cr},~j=1,2
\end{align}

The proof of the lemma is trivial. We omit here for brevity. Lemma 1 indicates that the linear approximation does not change $O_{11}$ and $O_{12}$ in the free flow state. 

The linearization of $O_{11}$ and $O_{12}$ can greatly accelerate the solution procedure of the proposed MPC model. 
Constraints Eq.(\ref{eq:n11d}) and Eq.(\ref{eq:muout}) can be reformulated as follows. 
\begin{align}
n_{11}(k+1) &= n_{11}(k) + D_{11}(k)C + \beta_{21}(k)C - \min\{ vCn_{11}(k), \gamma_{10}C + \gamma_{11}Cn_{11}(k)+\gamma_{12}Cn_{12}(k)\}, ~\forall 0\leq l\leq L-1   \label{eq:con1} \\
\mu_m^i(k) &= \min\Big\{\frac{\alpha_m^i}{|I|}vn_{12}(k), \frac{\alpha_m^i}{|I|}\Big(\gamma_{20}+\gamma_{21}n_{11}(k) + \gamma_{22}n_{12}(k) \Big), \sum_{p\in P^i}{s_{mp}^ig_p^i(k+l)}\Big\},~\forall i \in I, ~\forall 0\leq l\leq L \label{eq:con2}
\end{align}
By replacing Eq.(\ref{eq:n11d}) and Eq.(\ref{eq:muout}) with Eq.(\ref{eq:con1}) and Eq.(\ref{eq:con2}), the highly nonlinear and nonconvex MPC can be reformulated into a piecewise linear optimization problem. Though still non-convex, the problem can be solved using the existing algorithms (for example the branch and cut algorithm in \citet{keha2006}) by introducing binary variables which represent whether these constraints are chosen. Since there will not be too many binary variables, the problem can be solved in a reasonable time (30s for our simulated cases). 

We notice that the only constraint which still makes the problem non-convex is Eq.(\ref{eq:muother}) for $m\in M^i_{in}$, as the sign before $\beta_{21}(k+l|k)$ in Eq.(\ref{eq:con1}) is positive. To accelerate the procedure, we approximate $\mu_m(k+l|k) =\sum_{p\in P^i}{s_m^pg^p(k+l)}$, for $m \in M^i_{in}$. Two remarks should be mentioned for this approximation.  
\begin{enumerate}
\item[1)] This approximation might overestimate $\mu_m(k+l|k)$ if $x^i_m(k+l|k)/C + q^i_m(k+l|k)< \sum_{p\in P^i}{s_m^pg^p(k+l)}$, i.e. when demand to enter the network is low. However, in scenarios where there are always high demand, e.g. during peak hour, the approximated model is practically equivalent to the original model. Also note that by overestimating the network inflows, the obtained control policy can be more conservative at the network level. 
\item[2)] If the signal timing is very flexible (i.e. $g_{p,\min}^i$ is sufficiently small and each green phase discharges only one stream), $g^p(k+l)$ can be set as small as possible. In such cases, the approximated model is also equivalent to the original model. 
\end{enumerate}
With this approximation, we replace Eq.(11) with two constraints Eq.(\ref{eq:mua1}) and Eq.(\ref{eq:mua2}). 
\begin{align}
&\mu_m^i(k+l|k) = \min\{x^i_m(k+l|k)/C  + q^i_m(k+l|k), \sum_{p\in P^i}{s_{mp}^ig_p^i(k+l)}\},~\forall m \in M^i\backslash M^i_{\rm{out}}\backslash M^i_{\rm{in}}, ~\forall i \in I,~\forall 0\leq l\leq L-1 \label{eq:mua1}\\
&\mu_m^i(k+l|k) =\sum_{p\in P^i}{s_{mp}^ig_p^i(k+l)},~\forall m \in M^i_{\rm{in}}, ~\forall i \in I,~\forall 0\leq l\leq L-1\label{eq:mua2}
\end{align}
Eq.(12) is further replaced with Eq.(\ref{eq:12'}) to guarantee that $x^i_m(k+l|k)\geq 0$. 
\begin{align}
&x^i_m(k+l+1|k) = \max\{0,x^i_m(k+l|k) + q^i_m(k+l|k)C - \mu_m^i(k+l|k)C\},~\forall m \in M^i\backslash M^i_{\rm{out}},~\forall i \in I,~\forall 0\leq l\leq L-1 \label{eq:12'}
\end{align}

Having been transformed to its current form, the approximated model becomes a convex piecewise linear model. Both the traffic accumulations $n$ and the queue length $x$ can be represented as a sum of $\max$ functions of the green ratios $g$. This convex piecewise linear model can be equivalently reformulated into a linear programming problem by relaxing the $\min$ and $\max$ operator with the corresponding inequalities. The model can be solved efficiently with the classic simplex algorithm~\citep{nelder1965simplex} or an interior point algorithm~\citep{kojima1989primal} in a commercial solver, such as  Cplex. The average running time for solving the multi-scale problem in each cycle is 0.5s on a PC with one CPU core.

\section{Case study on the multi-scale controller}
\label{sec:4}
\subsection{Simulation settings}
\label{sec:4.1}
We simulate a typical morning-peak period. 
The MFD and the demand profile of this network are designed to mimic the aggregated traffic features of the city center of Zurich, Switzerland . 
The city has implemented a similar type of perimeter control~\citep{ortigosa2014study}. 

\begin{figure}[htbp]
\centering
\subfigure[Network level demand]{
\includegraphics[width=0.3\textwidth]{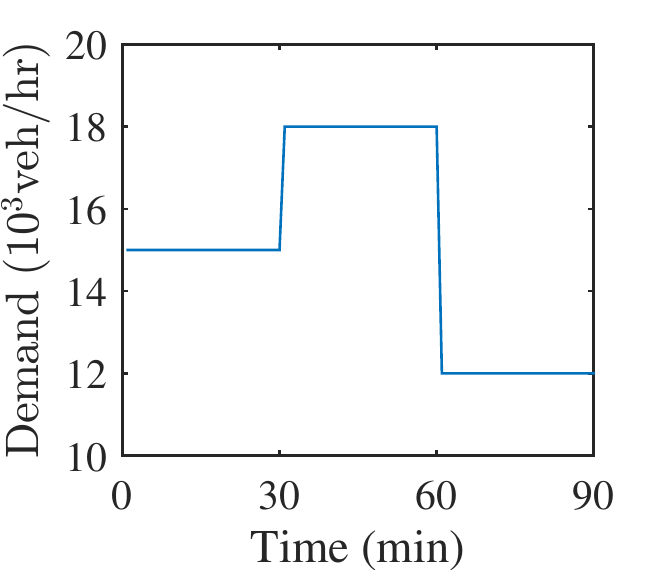}}~~
\subfigure[Intersection streams]{
\includegraphics[width=0.3\textwidth]{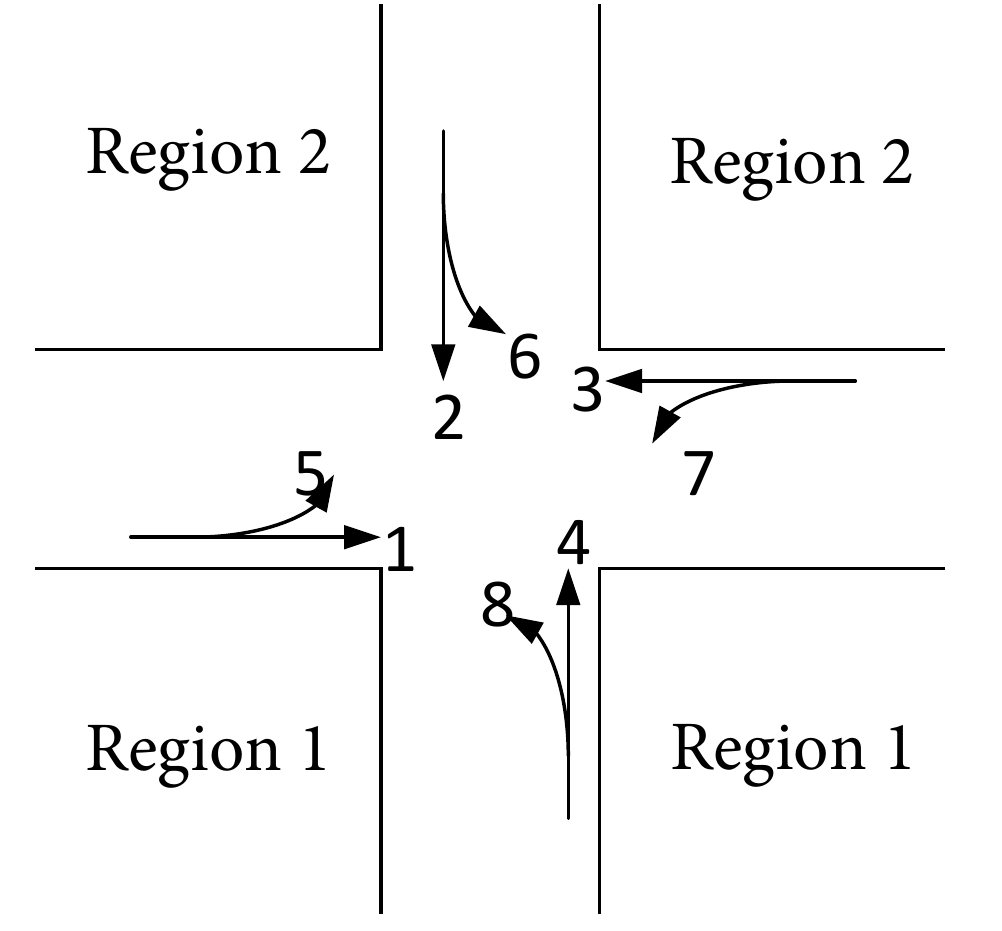}}~~
\subfigure[Intersection phases]{
\includegraphics[width=0.3\textwidth]{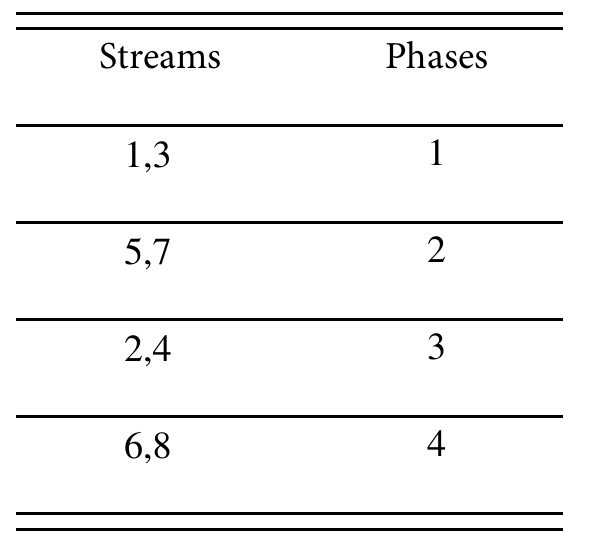}}
\centering
\caption{Simulations settings. At the intersection level (b and c), streams 2 and 7 are the inflow to the city network, and streams 4 and 8 are the outflow from the city network.}
\label{fig:settings}
\end{figure}
The total demand for entering the network during a 1.5 hour period is shown in Fig.\ref{fig:settings}a). Note that this demand profile generates a large amount of traffic accumulation for testing the controllers. 
The resulting MFD without the application of any type of perimeter control reaches oversaturated states (i.e. the network gets congested).

This network consists of 20 intersections at the perimeter. For the prediction model, we use a triangular MFD with $v=5\rm{hr}^{-1}, w = 2.5\rm{hr}^{-1}$ and $n_{cr}=3000$veh. 
The maximum flow of the network is 15000veh/hr. The initial traffic accumulation in the onset of the simulation is assumed to be 2000veh, below the critical accumulation.  Each intersection is assumed to have the phases and streams shown in Fig.\ref{fig:settings}b) and Fig.\ref{fig:settings}c)\footnote{Right turn is not included in this case study, as it will complicate our signal plan with limited added value. It is not difficult for our model to handle the cases with either protected or unprotected right turns.}. The cycle length $C$ is 1 min. The maximum allowed green ratio for each intersection $g_{\max}^i$ is 0.9, and the minimum green ratio $g_{p,\min}^i$for each phase is 0.1 (which means 6s). 
There is no offset between intersections.

To illustrate the performance of the proposed controller, we compare three control algorithms in the following case study: 1) the proposed multi-scale control algorithm, 2) the classical perimeter control algorithm that has been developed in the traffic control literature using Proportional, Integral and Differentiation (PID) approach \citep{aastrom2006advanced}, and 3) a bang-bang type of control algorithm which is also one of the most classical benchmarks in the control theory literature~\citep{bellman1956bang}. For both the PID controller and the bang-bang controller, we first determine the total inflow and outflow of the network. These flows are then distributed to each intersection based on the local queue information. A system with perfect information is considered in Section 4.2; Systems with moderate noises and large noises are considered in Section 4.3 and Section 4.4, respectively.  
The numerical simulation is conducted with Python, with integration of a Cplex C++ API for solving the model. The average solution time for one cycle is less than 1s.  
The moving horizon for the proposed multi-scale MPC is chosen as 20 cycles. 

\subsection{Performance of the multi-scale controller}
\label{sec:4.2}

\begin{figure}[t]
\centering
\subfigure[Multi-scale MPC]{\includegraphics[width=0.28\textwidth]{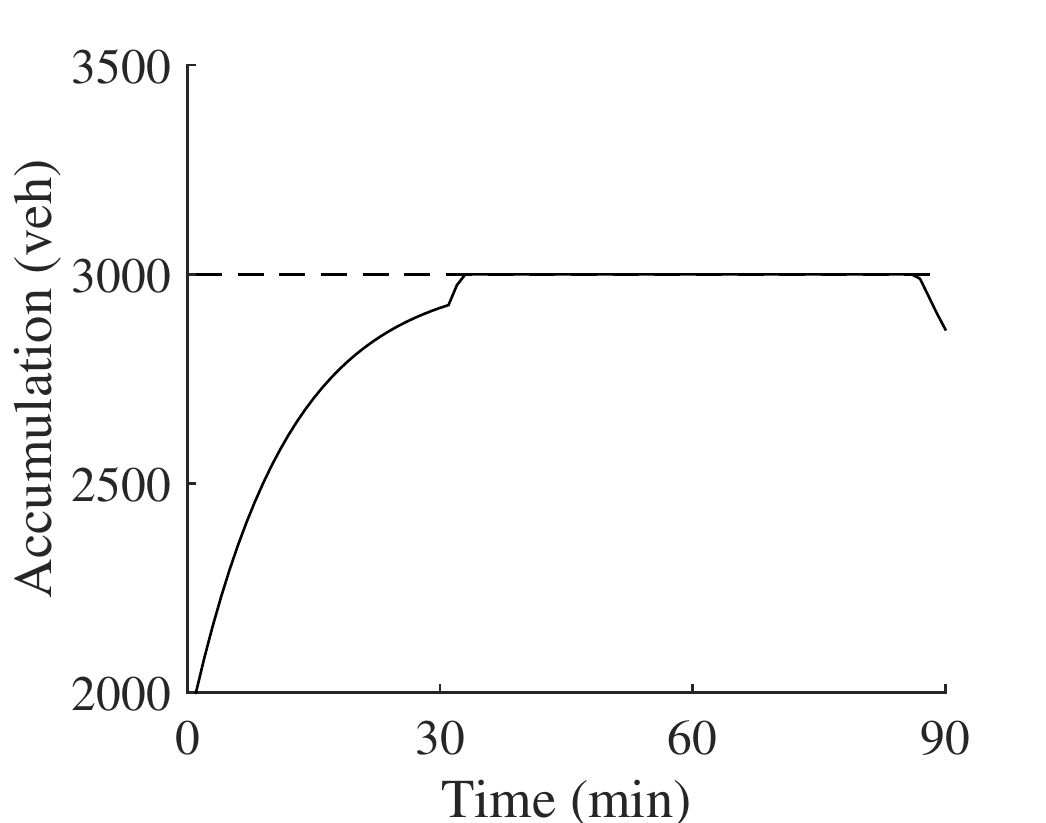}}
\subfigure[Classical PID controller]{\includegraphics[width=0.28\textwidth]{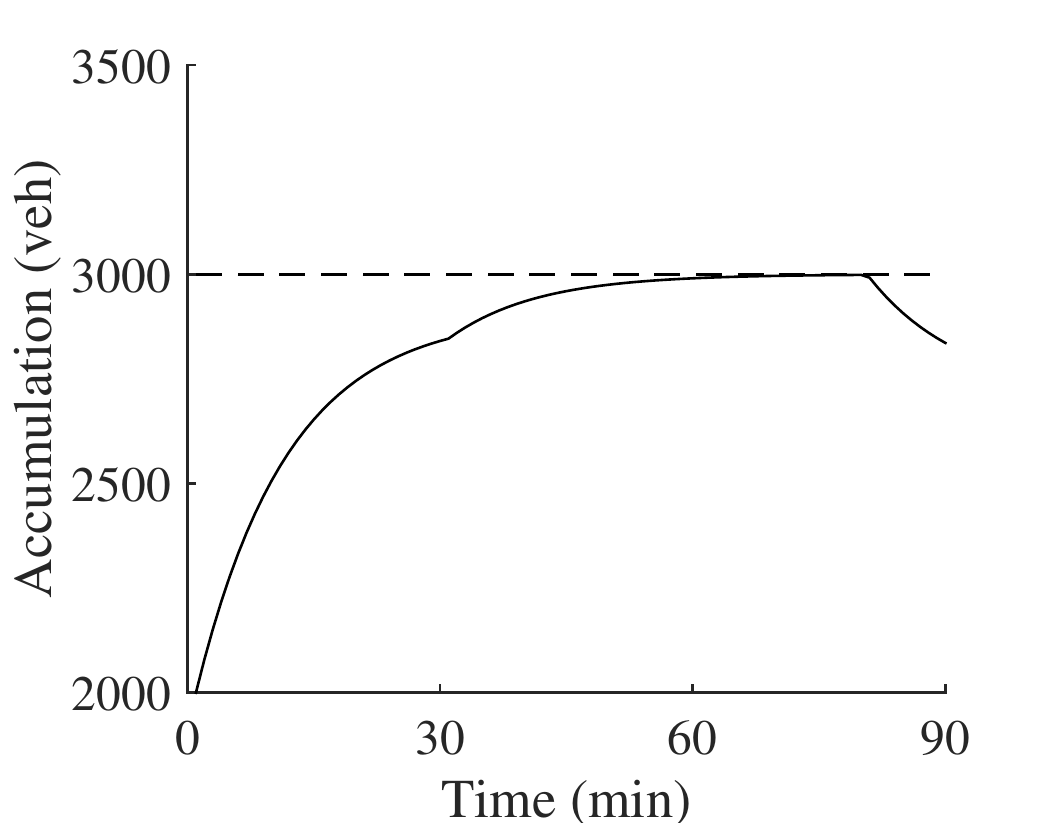}}
\subfigure[Bang-bang controller]{\includegraphics[width=0.28\textwidth]{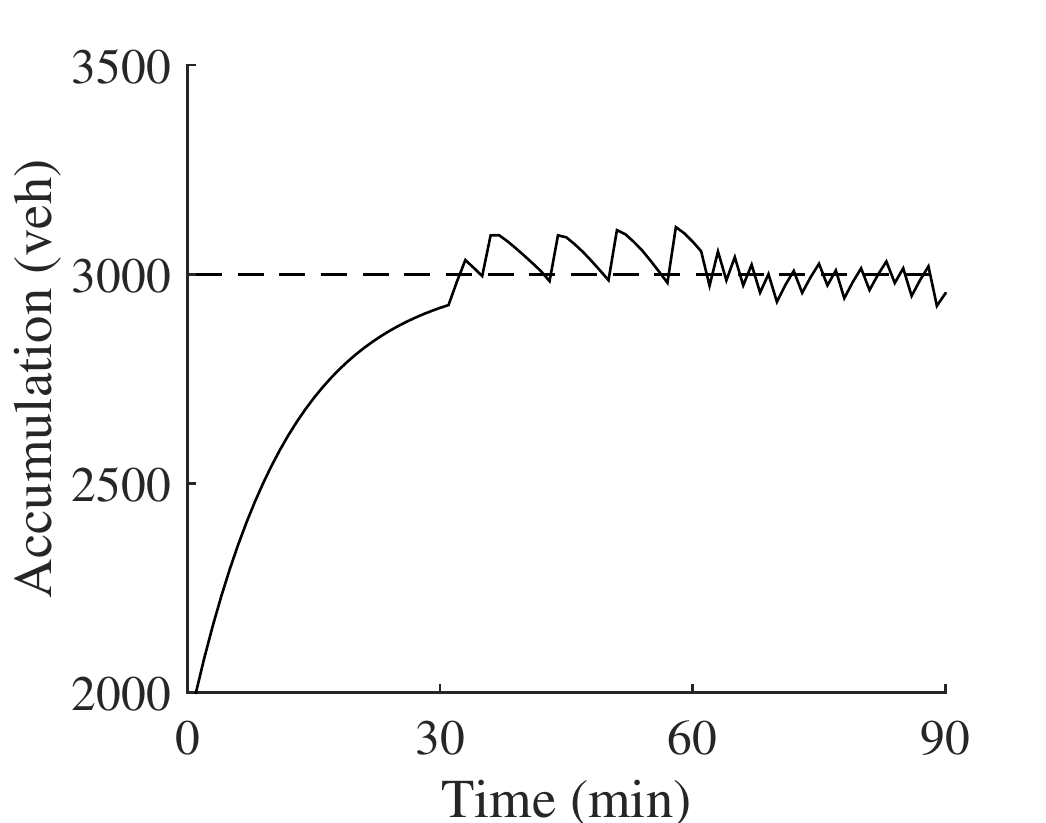}}
\centering
\caption{Accumulation comparison between the MPC approach, the classical PID approach and the bang-bang control in scenarios without noises.}
\label{fig:accumulation}
\end{figure}

\begin{figure}[htbp]
\centering
\subfigure[Onset of the peak hour (min. 35) at inflow streams]{
\includegraphics[width=0.4\textwidth]{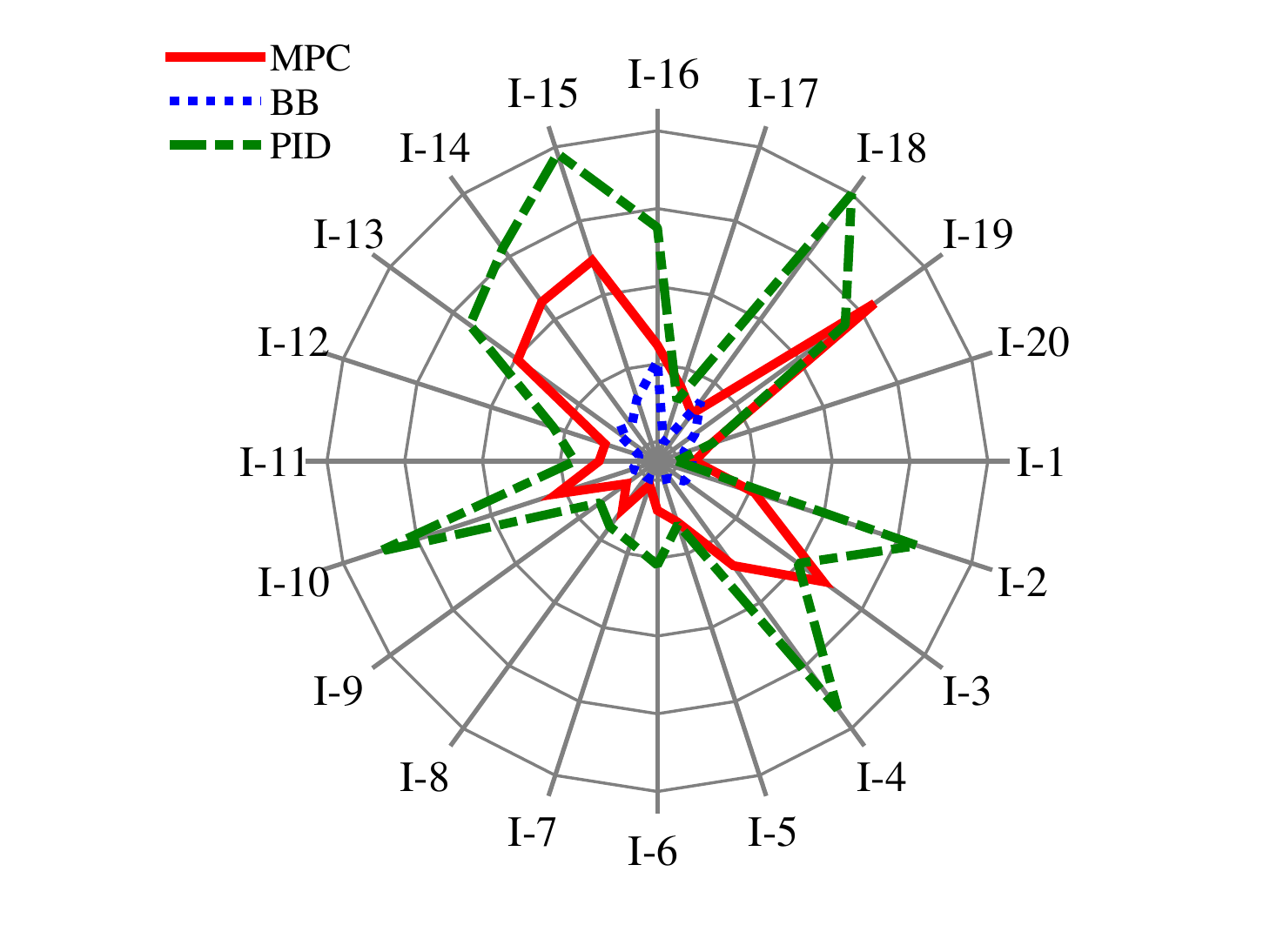}
}
\subfigure[Middle of the peak hour (min. 70) at inflow streams]{
\includegraphics[width=0.4\textwidth]{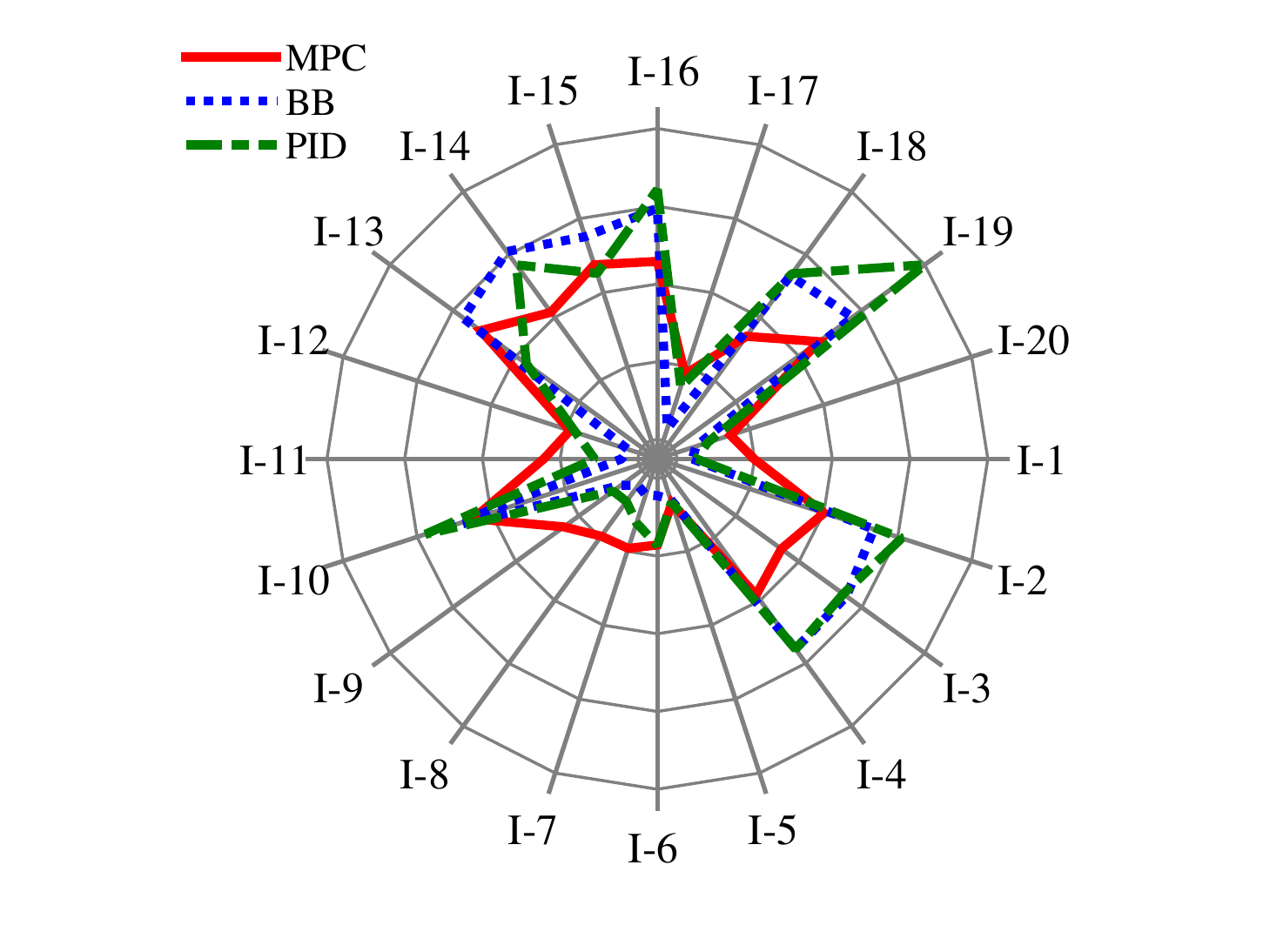}
}
\subfigure[Onset of the peak hour (min. 35) at side streams]{
\includegraphics[width=0.4\textwidth]{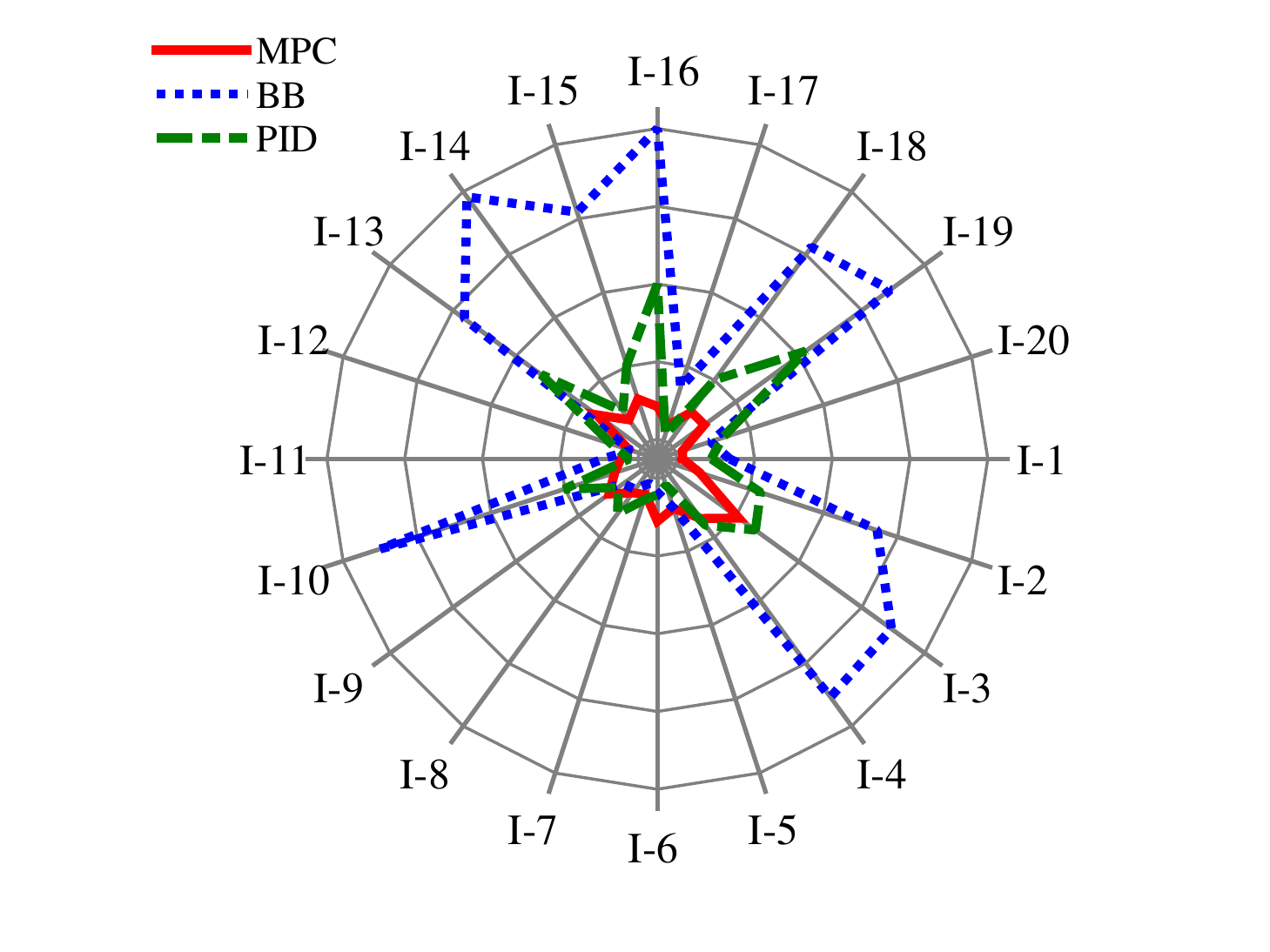}
}
\subfigure[Middle of the peak hour (min. 70) at side streams]{
\includegraphics[width=0.4\textwidth]{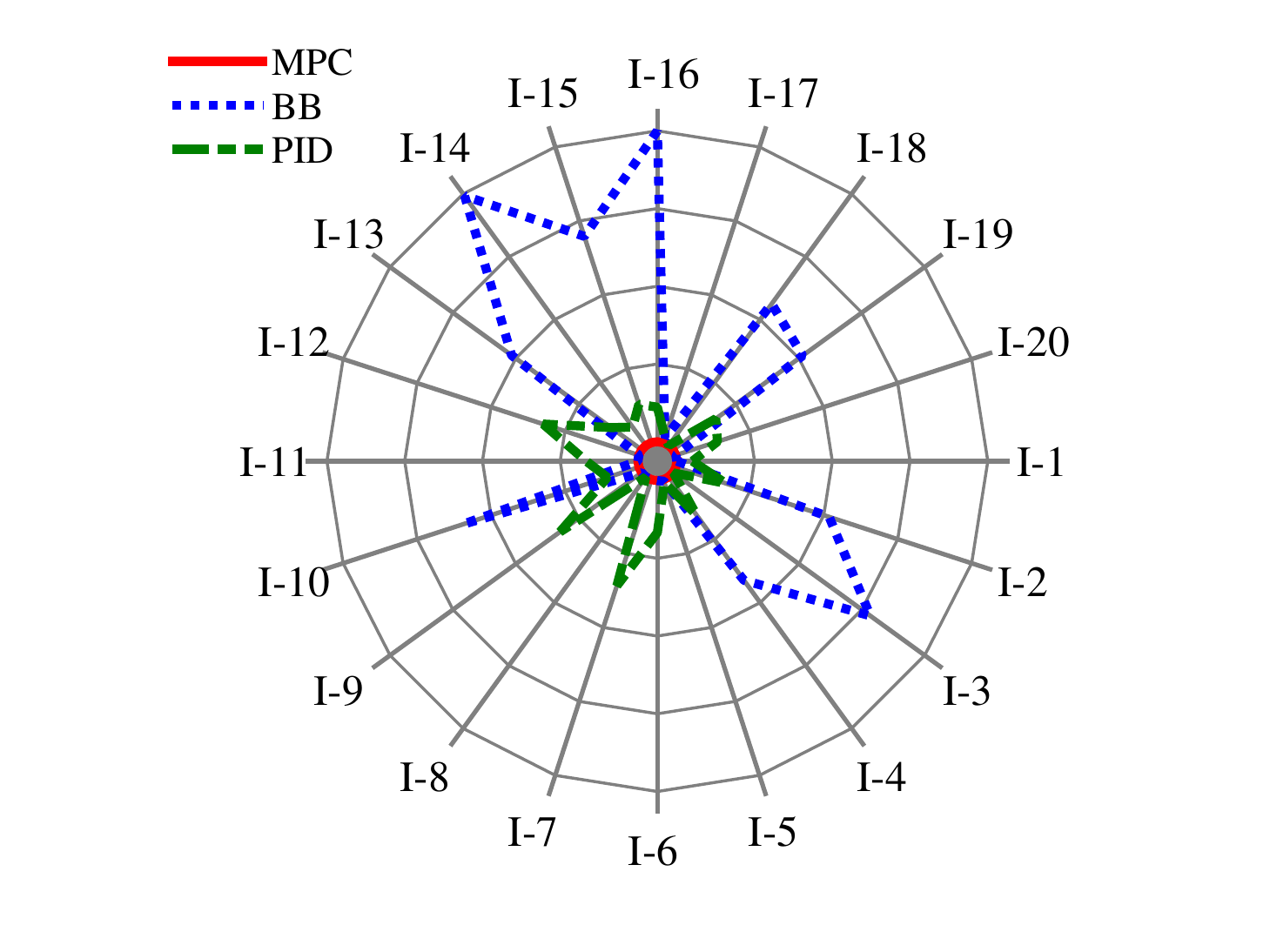}
}
\centering
\caption{Queue length comparison between the MPC approach, the bang-bang control approach, and the classical PID approach in scenarios without noises. Queue lengths are normalized with the minimum and maximum queue length in the network.}
\label{fig:queue_radar}
\end{figure}

\begin{figure}[htbp]
\centering
\includegraphics[width = 0.35\textwidth]{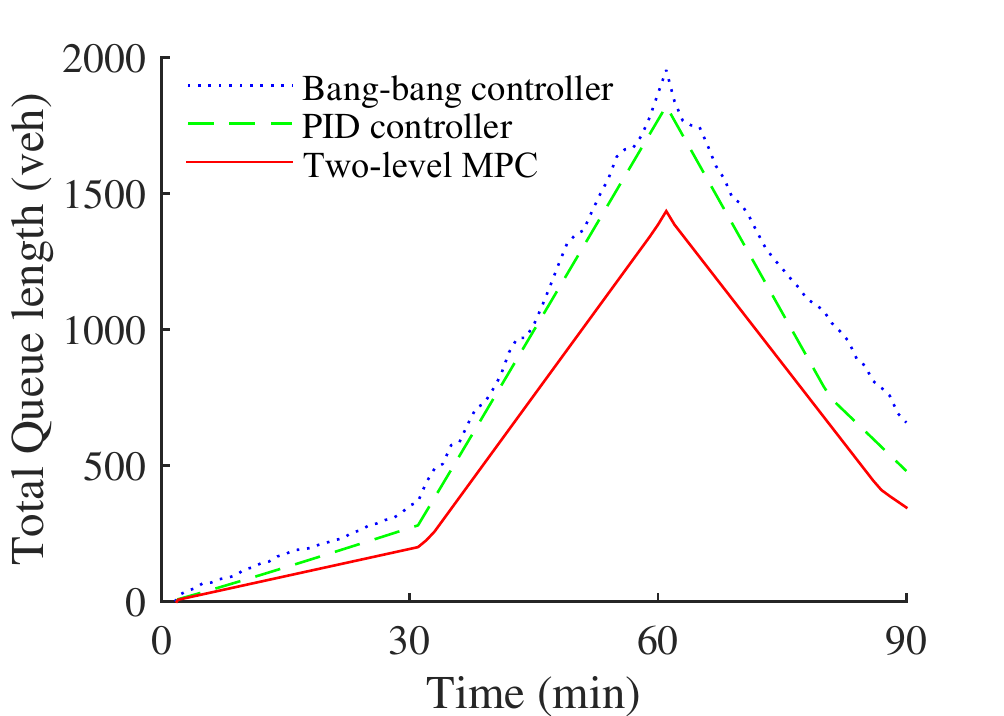}
\centering
\caption{Total queue length for all streams and intersections without noises.}
\label{fig:totalqueue}
\end{figure}

Let us define two performance measures, the traffic accumulation in the center and the queue length at the perimeter intersections. The queue values are normalized and dimensionless (rescaled in relation to the maximum and minimum queue values of the whole network), given that the intersections have different physical configurations (e.g. the discharging flow, arrival flow, etc). The resulting time series of the two measures are displayed in Fig. \ref{fig:accumulation}, Fig. \ref{fig:queue_radar} and Fig. \ref{fig:totalqueue} for the three controllers.

As shown in Fig.~\ref{fig:accumulation}, all three controllers succeed in effectively operating the city center network, by maintaining the accumulation around the critical value. 
It appears that the bang-bang control algorithm generates oscillations in the system dynamics. 
The classical bang-bang control algorithm follows an all-or-nothing rule, which is as expected to cause aggressive control actions towards the change of the system states.
The oscillations exist even in scenarios without noises in the control system (the type of noises considered later will be detailed in Section 4.3). 
The PID controller, regulates the system in a relatively less aggressive manner, and maintains an the accumulation slightly smaller than the critical accumulation. 
This is because the inflow and the outflow of the network, given by the PID controller, may not satisfy the constraints at the intersections, as the local level dynamics and signal settings are not explicitly considered.  
The proposed MPC controller appears to produce promising results, exhibiting high stability and efficiency at the network level. First, in comparison with the bang-bang controller, the MPC produces oscillation-free system dynamics. Second, in comparison with the PID controller, the MPC assures a higher network utilization, accommodating traffic at the critical accumulation level for a much longer period of time (roughly 30 minutes). Regarding the global performance, the three controllers seem to be equivalently functioning. After a careful analysis on the local level, however, the real advantage of the proposed controller is revealed. The differences are indeed found at the local intersections.

Fig. \ref{fig:queue_radar} displays radar-shape plots of queue length distribution  at the intersections. These plots illustrate the average queues in the inflow (streams 2 and 7) and side (streams 1, 3 ,5 and 6) direction to the network,  during two time periods: the beginning of the peak hour (Fig. \ref{fig:queue_radar}a and Fig. \ref{fig:queue_radar}c), and the middle of the peak hour (Fig. \ref{fig:queue_radar}b and Fig. \ref{fig:queue_radar}d). 
For the outflow direction, the demand is served by all of the three approaches, and there is almost no queue at any intersection.

The solid radial lines in gray represent the 20 intersections, while the rings are the queue length contour references. The radar pattern connects the length of the queue of the adjacent intersections, aiming at providing a clear illustration on the resulting queues under the different control schemes.   
Recall that the queue values have been normalized for comparable illustration.
As a general remark, the proposed multi-scale MPC controller handles the queue well for the concerned directions during both time periods. Comparing to the other two controllers, it generates shorter queues at the most intersections. 
From Fig.\ref{fig:queue_radar}a, it seems that the bang-bang controller works the best for managing the inflow to the network in the onset of the peak hour. 
This is not surprising, because the bang-bang controller gives the highest priority to the inflow streams by sacrificing the majority of traffic from the side streams (see Fig. \ref{fig:queue_radar}c). 
Moreover, the bang-bang controller may release queues of the inflow direction drastically which increases the network accumulation and negates the network-level performance (e.g. it can generate a traffic accumulation larger than the critical accumulation). 
In the middle 
 of the peak hour, it can be observed that the results are totally different. 
The MPC controller clearly works the best among the three controllers. As displayed in Figure \ref{fig:queue_radar}b , it reduces significantly the queues at most intersections. Furthermore it looks that the queues are more evenly distributed among intersections. 
The MPC controller outperforms the classical PID controller in both scenarios for the inflow queues.  
This is expected, as the PID controller does not take into account the detailed configuration and the queue of the intersections.
Furthermore, as observed in Fig.\ref{fig:queue_radar}c and Fig.\ref{fig:queue_radar}d, the proposed control improves the traffic operation of side directions， as well. 
It is found that the proposed controller significantly reduces the local queues, by up to 60\% compared to the bang-bang and the PID approaches. 
The bang-bang controller, not surprisingly, can result in much longer queues, e.g. for intersection 16 (I-16) the queue length is double under bang-bang as the one under the MPC or the PID.

Fig.~\ref{fig:totalqueue} shows the total queue length at each cycle for the three algorithms. It is evident that the proposed controller renders smaller queues than the other two algorithms.

We have also carried out analysis for scenarios where the demand is low. Under such conditions, the performance is similar for all controllers, e.g. the difference in the resultant total delays is ignorable. This is reasonable, since the magnitude of the impact of local queue on global performance is directly determined by the level of demand flowing into the network.

\subsection{System performance under moderate noises}
\label{sec:4.3}

In this section, we test the performance of the proposed approach in a noisy environment. 
Precisely, we consider the measurement noises, prediction uncertainties, and stochasticity in the MFD\footnote{ Note that the noisy environment also includes scenarios with uncontrolled intersections. The noises brought by the uncontrolled intersections \citep[as in][]{Keyvan-Ekbatani2012} can be combined with the internal demand.  }. 
Recall that we assume the connected vehicles are the only information source. 
In a network fully equipped with connected vehicles, perfect information can be provided at both network and local levels for the controller. 
However, when the penetration rate is insufficient, an inaccurate state measurement (i.e. traffic accumulation and queue lengths) would be expected. 
These noises affect the accuracy of the initial state in the MPC approach (see Fig.\ref{fig:robust}), i.e. $\hat{n}_{11}(k),~\hat{n}_{12}(k),~\hat{x}^i_m(k)$ in the proposed model Eq.(\ref{eq:obj})-Eq.(\ref{eq:xinit}), thus control efficiency could be hindered. 
Furthermore, as the traffic  dynamics can be highly stochastic, it is difficult to predict the future arrivals and departures with complete accuracy. 
 
It is then important to evaluate how less than perfect information regarding these variables affects the control approach.

\begin{figure}[htbp]
\centering
\subfigure[Green ratios at a representative intersection]{
\includegraphics[width=1\textwidth]{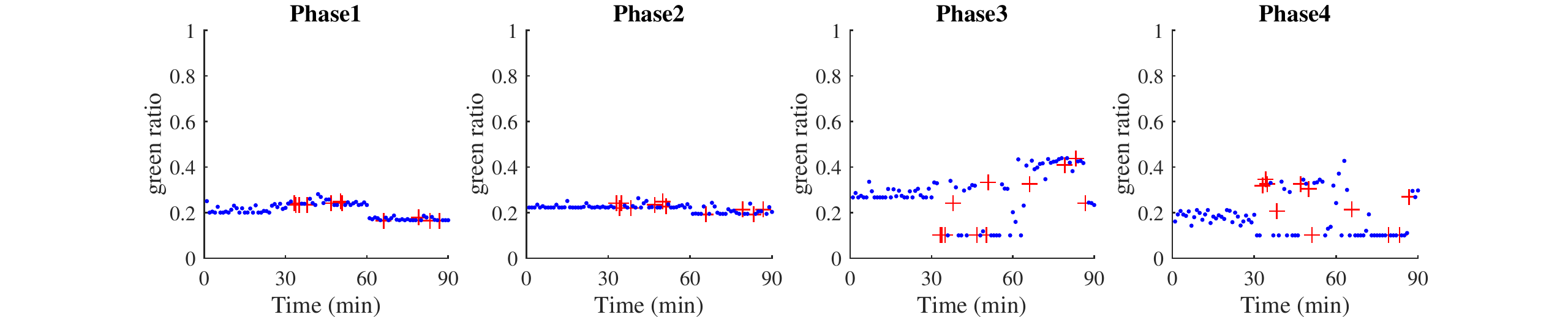}}
\subfigure[Departure flows (across all intersections)]{
\includegraphics[width=0.75\textwidth]{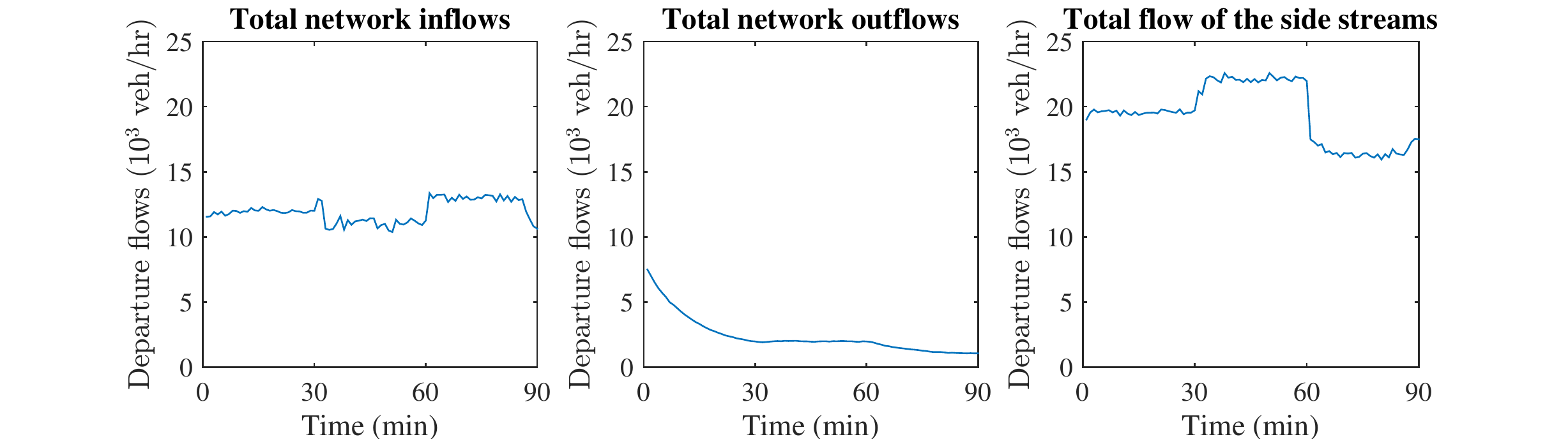}}
\subfigure[Network traffic accumulation]{
\includegraphics[width=0.25\textwidth]{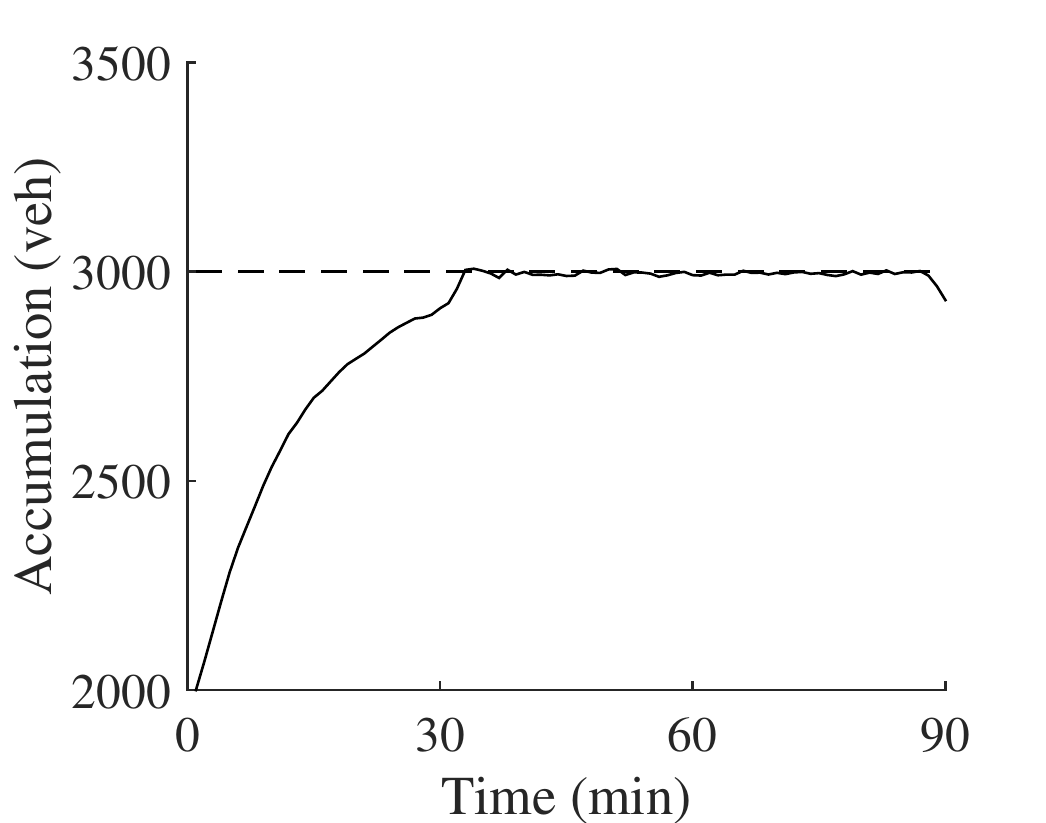}}
\centering
\caption{System performance with moderate noises}
\label{fig:moderate}
\end{figure}

In this section, we look at the system performance for moderate noise level. 
The error for the traffic accumulation and queue length measurement is assumed to follow a Gaussian distribution with a mean of 0 and a standard deviation of 5\% of the mean value.   
Furthermore, we simultaneously consider the stochasticity existing on the MFD, so the outflow from the network has e.g. a 10\% fluctuation (i.e, the actual outflow from the network follows a uniform distribution between $0.9O(t)$ and $1.1O(t)$). The errors in the demand prediction are assumed to follow a Gaussian distribution with mean 0 and a standard deviation of 10\% of the value. 
These three noise variables are assumed to be independent.



Fig.\ref{fig:moderate}a displays the resulting green ratio over time, for the four different phases of a representative intersection operating under the perimeter control algorithm. 
The performance of the other intersections looks similar even though they have different demand patterns and capacities.  
Recall that Phases 2 and 4 accommodate the network inflow and outflow respectively, while Phase 3 serves both directions. 
Phase 1 regulates the directions other than the network inflow and outflow. It receives green time allocation up to 25\% of the cycle time. 
Such result is obtained because the controller aims to minimize the total delay of all traffic. 
A simplification of the intersections with only inflow and outflow direction cannot reflect the impact of control on local queues, and possibly overestimates the transfer flows between the two regions.
Treating the entire intersection is thus important. 
From the time series of the Phase 2, it can be observed that the restriction of the inflow to the network occurs during the peak hour between minutes 30 and 90.

Fig. \ref{fig:moderate}b shows the controlled inflow (Streams 2 and 7), outflow of the center region (Streams 4 and 8), and other controlled departure flows (Streams 1, 3, 5 and 6) at the aggregated level for the whole network. 
It is evident that all the departure flows are quite smooth under moderate noises. 
Fig. \ref{fig:moderate}c displays the time series
of traffic accumulation. The changes of flow restriction and traffic accumulation at the network level appear to be stable and smooth.

\subsection{Performance deterioration  with strong noises in the system}
\label{sec:4.4}

We investigate now the performance when the system error increases. The corresponding noises on the measurement of accumulation and queue length reach to 15\% (instead of 5\% in Section \ref{sec:4.3}). 
The errors in the demand prediction are assumed to follow a Gaussian distribution with mean 0 and a standard deviation of 30\% (in contrast to 10\% in Section \ref{sec:4.3}) of the value. 
The outflow calculated from the MFD is assumed to have a 20\% fluctuation (in contrast to 10\% in Section \ref{sec:4.3}) .
All the other settings remain unchanged.

\begin{figure}[htbp]
\centering
\subfigure[Green ratios at a representative intersection]{
\includegraphics[width=1\textwidth]{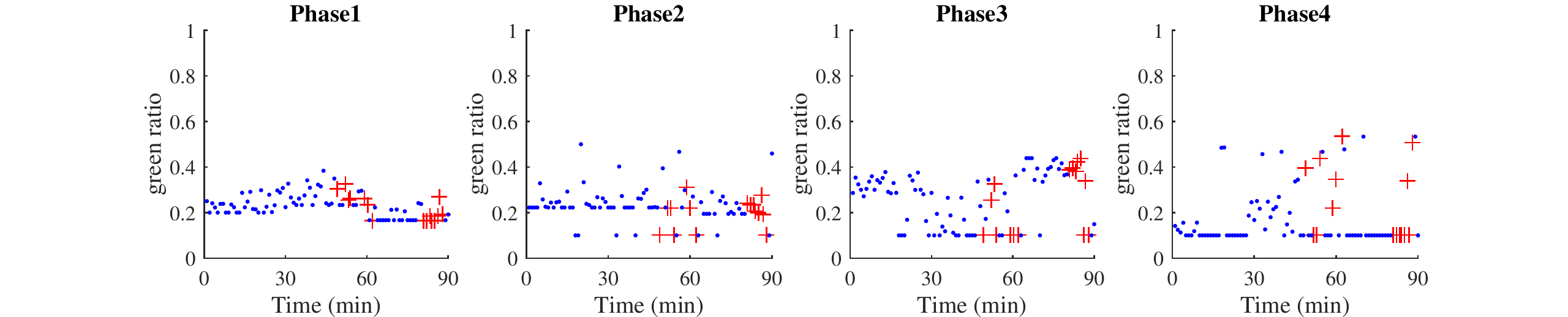}}
\subfigure[Departure flows (across all the intersections)]{
\includegraphics[width=0.75\textwidth]{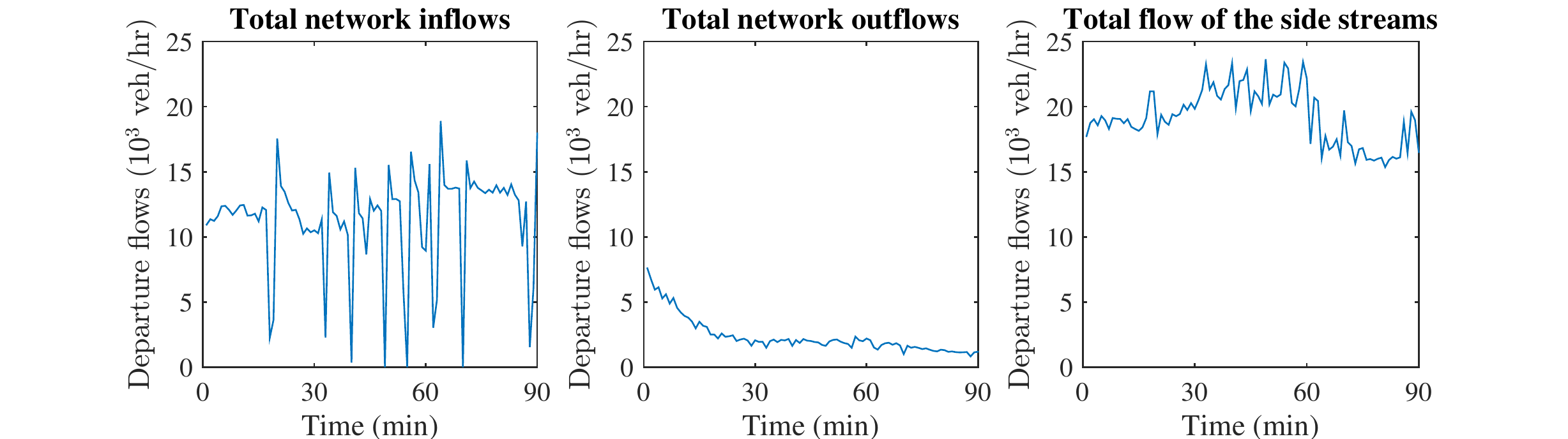}}
\subfigure[Network traffic accumulation]{
\includegraphics[width=0.25\textwidth]{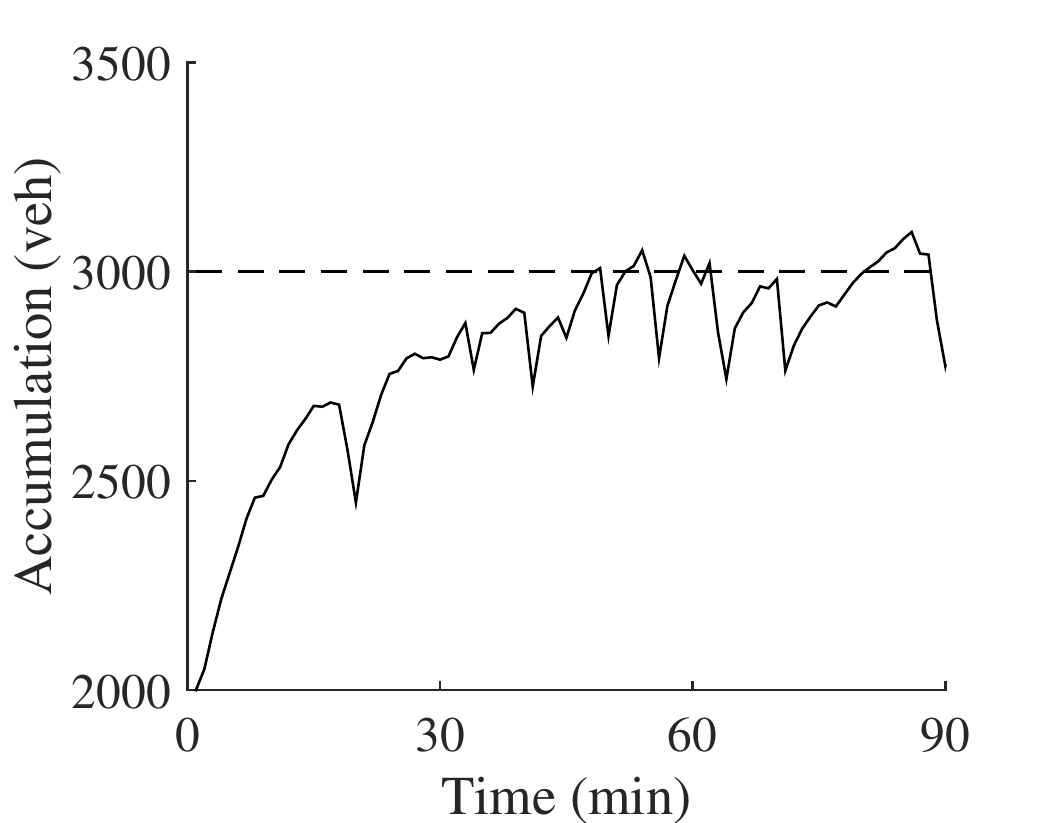}}
\centering
\caption{System performance with strong noises.}
\label{fig:large}
\end{figure}

We reproduce the same graphs from Fig. \ref{fig:moderate} in Fig. \ref{fig:large}. 
The accumulation time series in Fig.\ref{fig:large}c) clearly indicates a strong oscillation after the control is activated. 
Oscillations are also observed in the controlled flows. 
After a first glance, it may seem the proposed controller works in a bang-bang style, as up-and-downs of departure flow can be found at local intersections and at the network.  

To investigate in detail the difference in  mechanism between the proposed controller and the bang-bang controller, we plot in Fig.\ref{fig:large}a) the time series of the green ratio with scatters and colors. Red color represents the time intervals when the accumulation is larger than the critical accumulation. It can be seen that the proposed control does not follow a "green-or-red" logic which is employed by the bang-bang controller. The proposed controller sometimes even restricts the inflow, when the accumulation is below the critical one (for example at minute 44). As the proposed controller predicts the possible increase in that may cause queues later, inflow is restricted in advance. This reflects the fundamental difference of the proposed control from the bang-bang which takes only myopic decisions.

\begin{table}[htbp]
\renewcommand\thetable{2}
\centering
\caption{Comparison between the two MPC, the PID controller and the Bang-bang controller. Numbers within parenthesis reflect the change with respect to the multi-scale MPC for each of the scenarios. }
\label{my-label}
\begin{tabular}{l|l|l|l|l}
\hline                                                                                           &                  & Large noise & Moderate noise & No noise    \\ \hline
\multirow{3}{*}{\begin{tabular}[c]{@{}l@{}}Total\\   travel cost (hr)\end{tabular}}        & Multi-scale MPC & 7263 & 6122 & 5443 \\ \cline{2-5} 
                                                                                           & PID    controller           & 7722 (+6.3\%)
& 6663 (+8.8\%) & 5708 (+4.9\%)\\ \cline{2-5} 
                                                                                           & Bang-bang controller         & 8928 (+22.9\%)& 7481 (+22.2\%)& 6211 (+14.1\%)\\ \hline
\multirow{3}{*}{\begin{tabular}[c]{@{}l@{}}Total intersection\\   delay (hr)\end{tabular}} & Multi-scale MPC & 2536 & 1495 & 1082 \\ \cline{2-5} 
                                                                                           & PID controller              & 2950 (+16.3\%)& 1972 (+31.9\%)& 1341 (+23.9\%)\\ \cline{2-5} 
                                                                                           & Bang-bang controller         & 3953 (+55.9\%)& 2624 (+75.5\%)& 1485 (+37.2\%)\\ \hline
\end{tabular}
\label{tab:comparison}
\end{table} 

A quantitative comparison between the proposed multi-scale MPC, the PID controller, and the bang-bang controller in scenarios with large noises, moderate noises and no noises is shown in Table \ref{tab:comparison}. It is shown that the hybrid controller outperforms both the PID and the bang-bang controller, especially with noises. Specifically, the savings in network travel times is more than 4.9\% compared to the PID controller and more than 14\% compared to the Bang-bang controller. The savings in intersection delay is more than 16\% compared to the PID controller and more than 37\% compared to the Bang-bang controller.  

In summary, we have demonstrated that although the algorithm performs well under moderate noises, strong noises and uncertainties in the system due to the low penetration rates of CVs can create strong fluctuations in system dynamics.
This type of fluctuation indicates an over-correction of the controller, which hinders the performance of the system. 
To address this issue, we propose a stochastic control approach in the next section.

\section{Stochastic Controller}
\label{sec:5}
\subsection{Design of the stochastic MPC }
\label{sec:5.1}

Given the results of Section \ref{sec:4.4}, it is evident that the controlled performance of the system can become worse when perfect and accurate information on the traffic accumulations and demand is not available. 
Furthermore, utilizing deterministic traffic models (e.g. the non-scattered MFD) overlooks the stochastic nature of traffic dynamics, thus might make the control actions given by Eq.(\ref{eq:obj})-Eq.(\ref{eq:xinit}) non-optimal.

To this end, we extend the developed controller to account for the noises. 
We explicitly consider two types of noises: demand prediction errors and measurement errors in traffic accumulation and queue length. 
In other words, we consider the predicted demand ($D_{ab}(k+l|k)$, $q_m^i(k+l|k)$), and measurement ($\hat{n}_{ab}(k),~\hat{x}_m^i(k)$) to be random variables following known distributions. 
In fact, the distributions of these variables can be estimated from real-time or historical CV data with Bayesian filters (e.g. Kalman filter) or machine learning techniques. The details are beyond the scope of this paper. Interested readers can refer to \citet{yuan2012real}, \citet{gayah2013using}, or \citet{ramezani2015queue}. 

In contrast with existing robust controllers which optimize the worst-case performance of the controller in a noisy environment~\citep{haddad2014robust,haddad2015robust}, we employ a stochastic MPC based approach. For the definition of stochastic MPC, interested readers can refer to~\citet{kouvaritakis2004recent} and \cite{couchman2006stochastic}.
The reasons for proposing a stochastic MPC in this case are two-fold. 
First, the robust control scheme normally assumes bounded noises, which does not necessarily apply for this system. The actual internal demand at the network level and the arrival flows at the intersection level can be in a large range, which could make the resulting control action too conservative due to the strong stochasticity in the system. Second, applying the robust control scheme to this problem could render a complex, or even intractable optimization problem, considering the scale and complexity of the formulated model. 

Therefore, we propose a two-stage stochastic MPC approach, which employs stochastic programming. 
This type of control method has been extensively applied in many disciplines \citep[e.g.][]{Farina201653, Parisio201624,Tong2015474}. 
Recall that in the MPC approach, only the control action (green ratios) in the first cycle is executed (the others are for prediction purposes). 
Hence, we aim to determine the green ratios in the first optimized cycle (i.e. the current cycle) such that the optimal green ratios minimize the average cost in the $L$ future cycles. 
In the rest of this section, we establish the two stages of the stochastic MPC. 
The first stage master problem considers the green ratios to be executed (i.e. $g_p^i(k)$) subject to the physical constraints of the green ratios. The second stage subproblem estimates the total travel cost resulting from the chosen green ratios for each scenario arrival flow, internal demand, traffic accumulation, and  queue lengths for each sample in the sample space $R$. 

The first stage master problem is formulated as Eq.(\ref{eq:obj_1})-Eq.(\ref{eq:con_sto_1}). 
\begin{align}
\min~&E_rJ_{D,r} \label{eq:obj_1}\\
\mbox{s.t.~}& \sum_{p\in P}g_p^i(k) \leq g_{\max}^i, ~\forall i \in I \\
& g_p^i(k) \geq g_{\min,p}^i,~\forall
i \in I, p \in P^i \label{eq:con_sto_1}
\end{align}
where $r$ is one sample in the sample space $R$. $J_{D,r}$ is the optimal objective value of the second stage problem resulting from the decision variable $g_p^i(k)$ in sample  $r$. 

The second stage subproblem is formulated as follows.
\begin{align}
\min~J_{D,r} = C\sum_{l=1}^L{\Big(n_{11,r}(k+l|k)+n_{12,r}(k+l|k)+\sum_{i\in I}{\sum_{m \in M^i\backslash M^i_{\rm{out}}}x_{m,r}^i(k+l|k)}} \Big)
\end{align}
where the decision variables of the second stage problem in each sample $r$ are $g_{p,r}^i(k+l), 1\leq l\leq L$. The corresponding $n_{11,r}(k+l|k),~n_{12,r}(k+l|k),~x_{m,r}^i(k+l|k)$, together with $\mu_{m,r}(k+l|k)$, $\beta_{11,r}(k+l|k)$ and $ \beta_{12,r}(k+l|k)$, are regarded as functions of $g_{p,r}^i(k+l)$. These variables satisfy constraints Eq.(\ref{eq:n11d})-Eq.(\ref{eq:xinit}). 
Note that for the following cycles (starting from $k+1, l\geq 1$), the green ratios $g_{p,r}^i(k+l)$ are different for each sample $r$. In contrast, the green ratios in the current cycle, $g_{p}^i(k)$, are given as known variables for the second stage problem, and they are the same for each sample, i.e.
\begin{align}
g_{p,r}^i(k) = \hat{g}_p^i(k), ~\forall p \in P^i, \forall i \in I,\forall 0\leq l\leq L-1\label{eq:ginitr} 
\end{align}
We can simplify the second stage problem into a linear programming problem with the same approximation as utilized in Section 3.2.

The two stages interact with each other. The connection between the first stage problem and the second stage problem is two-fold. First, the objective function of the first stage problem is calculated by all the second stage problems. Second, the solution to the first stage problem is the initial solution for the second stage problem. 

However, as the sample space $R$ is continuous, it is not possible to enumerate all the samples to calculate the expected travel cost. One general approach to handle this problem is the so-called Sample Path Optimization~\citep{robinson1996analysis}, i.e. to sample a finite number of scenarios, $R_0$, following the given distribution. 
Specifically, we sample the arrival flow, internal demand, traffic accumulation, and  queue lengths as $q_{m,r}^i(k+l|k)$, $D_{ab,r}(k+l|k)$, $\hat{n}_{ab,r}(k)$ and $\hat{x}_{m,r}^i(k)$. Ideally, we should also adopt a large sample size (i.e. number of scenarios with different realizations of the same noise level) to estimate the expected travel cost more accurately. 
However, for the sake of computational efficiency, we limit the sample size to 20-50, as this number of samples is enough to cover the sample space. For information on how to generate representative samples, interested readers can refer to \citet{SOBOL196786} and \citet{ge2014efficient}. It is shown in Section 5.2 and Section 5.3 that the stochastic MPC performs sufficiently well with a sample size of 20 for the cases with similar city size to Zurich. Empirical simulations also indicate that the benefit of having a larger sample size than 20 is marginal.

We combine the two stages into one linear programming problem and solve it with the simplex method in Cplex. This is proved to be more efficient than the $L-$shape algorithm~\citep{higle1991stochastic} which iteratively solves the two stage problem.

 \subsection{Performance of the stochastic MPC}
 \label{sec:5.2}

\begin{figure}[t]
\centering
\subfigure[Green ratios at a representative intersection]{
\includegraphics[width=1\textwidth]{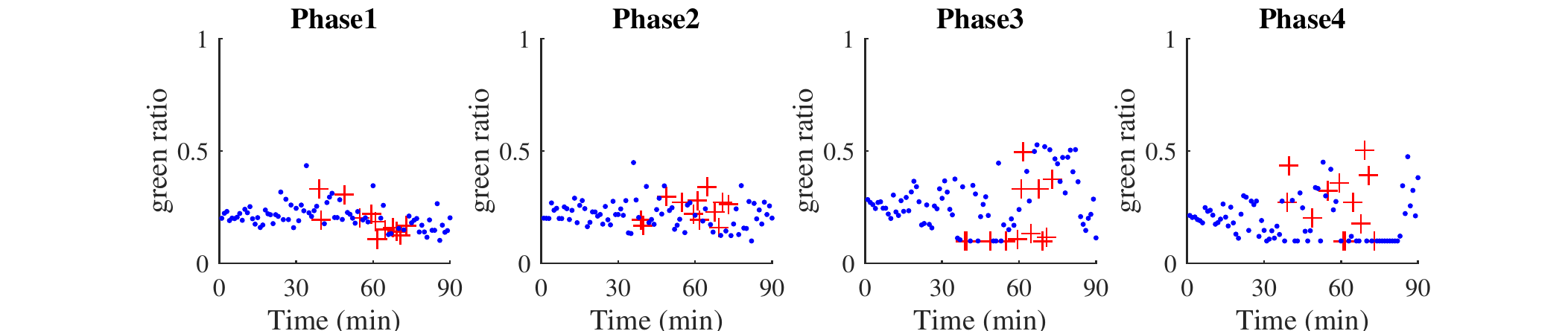}}
\subfigure[Departure flows (across all the intersections)]{
\includegraphics[width=0.75\textwidth]{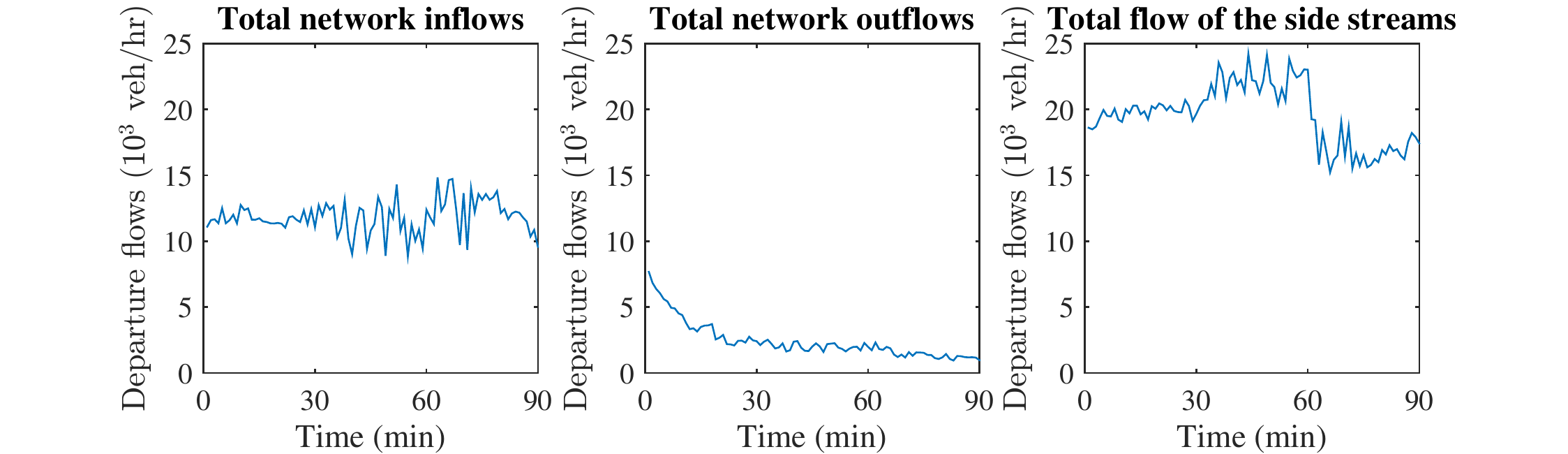}}
\subfigure[Network traffic accumulation]{
\includegraphics[width=0.25\textwidth]{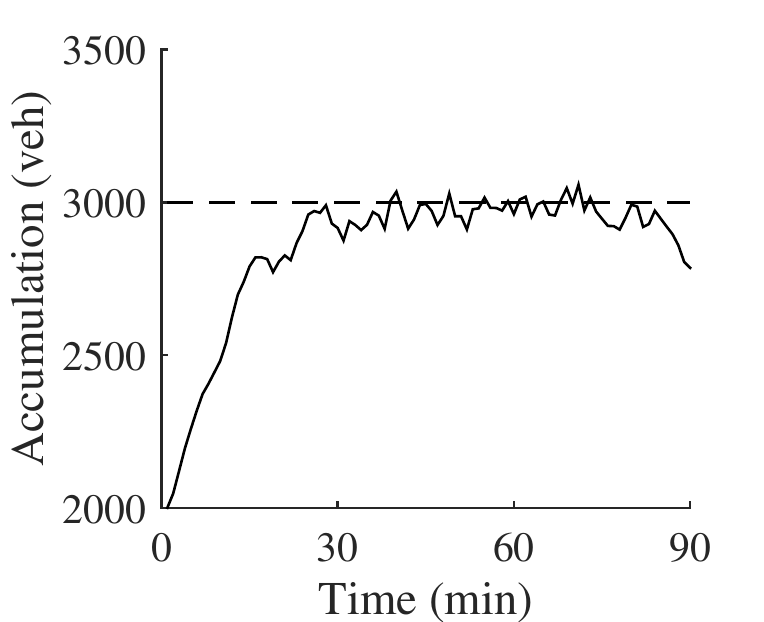}}
\centering
\caption{System performance of the stochastic controller with strong noises}
\label{fig:robust_result}
\end{figure}
We implement the stochastic controller using the same simulation settings (with strong noises) as in Section \ref{sec:4.4} and evaluate the performance of the controller. 
Fig. \ref{fig:robust_result} shows the comparison of the accumulation time series after applying the stochastic controller. 
The sample size is chosen as 20 and the computational time for each time step is 15s (in contrast with less than 1s for the multi-scale MPC based controller).  
Compared to Fig.\ref{fig:large}, it is evident that  the system changes in a more stable manner thanks to the stochastic controller. 
Compared to the oscillations of the multi-scale controller, the stochastic controller successfully maintains the accumulation around the critical one (Fig.\ref{fig:robust_result}b). The fluctuation of the network inflows is also reduced (Fig.\ref{fig:robust_result}c).

\subsection{Value of connected vehicles}
\begin{figure}[t]
\centering
\subfigure[Total travel cost]{\includegraphics[width=0.25\textwidth]{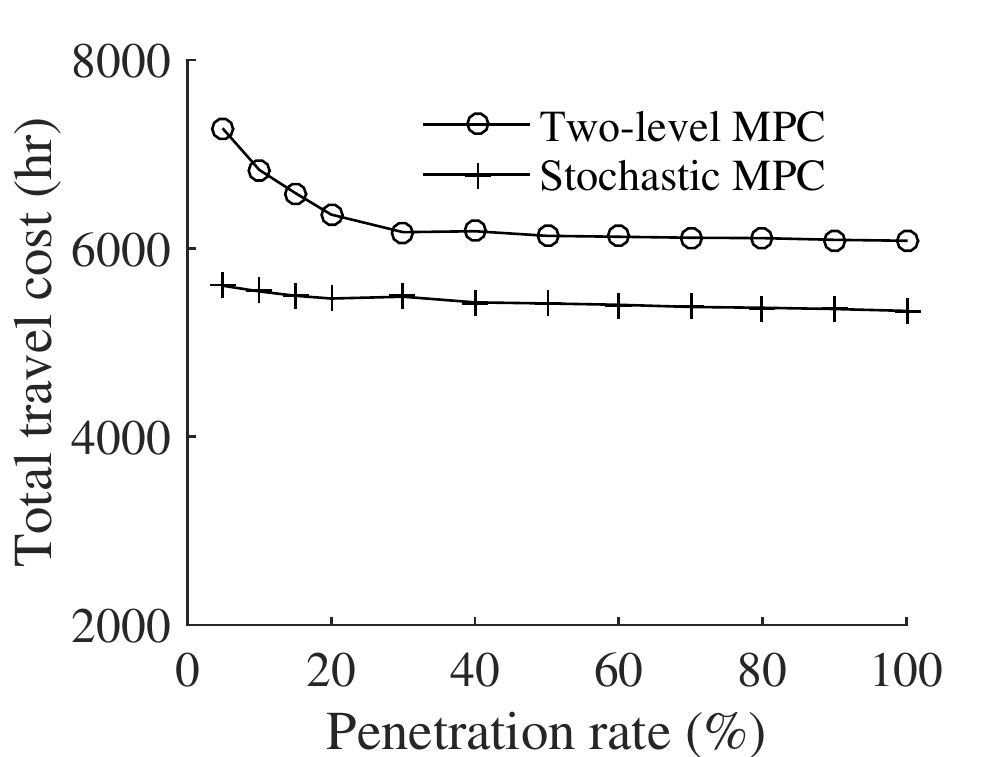}}
\subfigure[Network travel time]{\includegraphics[width=0.25\textwidth]{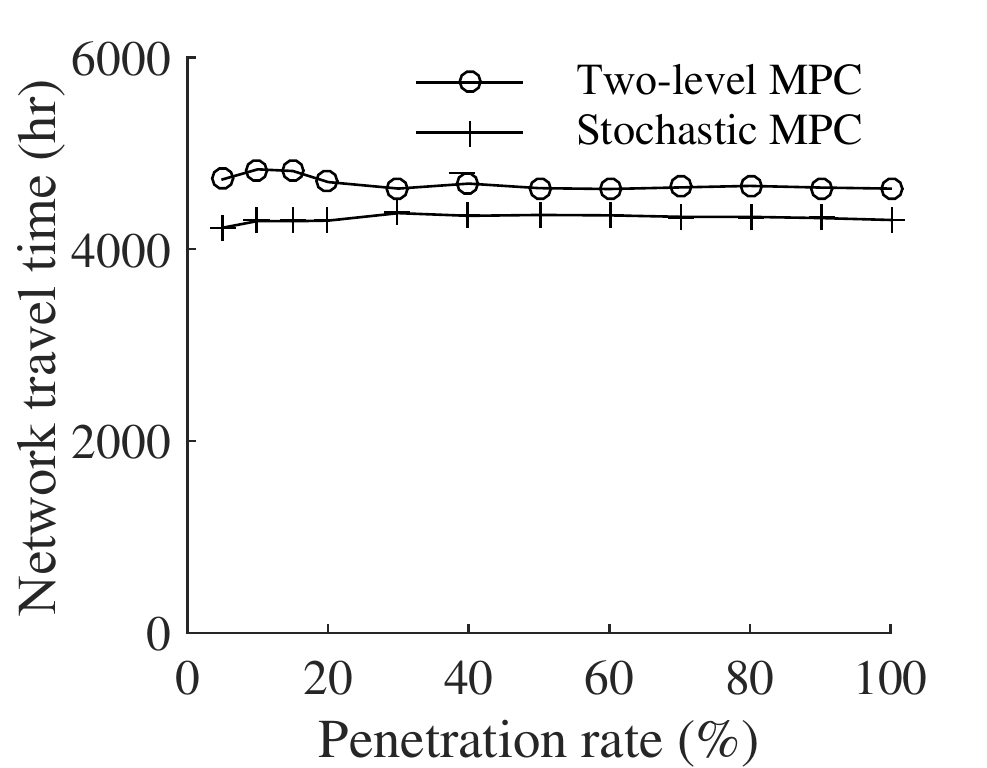}}\subfigure[Total intersection delay]{\includegraphics[width=0.25\textwidth]{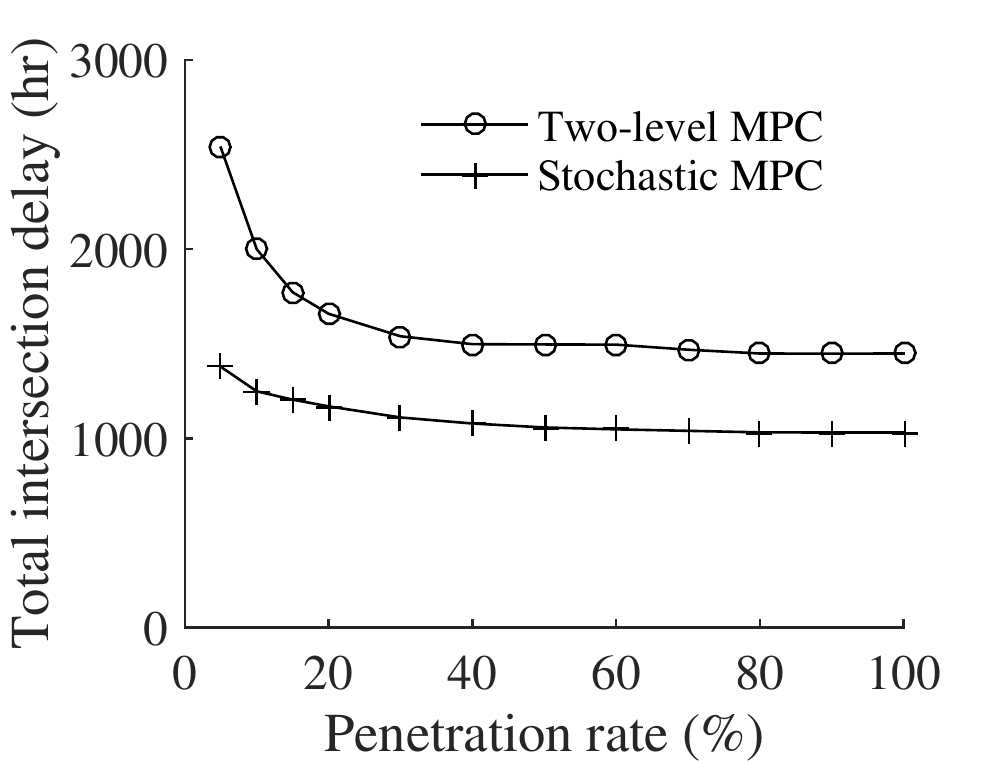}}
\caption{Performance of the stochastic MPC in different penetration rates of connected vehicles.}
\label{fig:cv}
\end{figure}

We evaluated the value of connected vehicles in this subsection. The penetration rate of connected vehicles affects the measurement of the traffic accumulation and the queue lengths. 
Existing works identified the relationship between the measurement errors and the penetration rates, e.g. in \citet{ramezani2015queue} and 
\citet{gayah2013using}. To avoid randomness, we select 10 random seeds for executing the controllers, we then compute the average value of each available performance measure.
The resulting performance measures from both the multi-scale controller and the stochastic controller are illustrated in Fig.\ref{fig:cv}.

Fig.\ref{fig:cv}a) shows the total travel costs
in the whole system over different levels of penetration rates. It is observed that the stochastic MPC always outperforms the multi-scale MPC for all the penetration rates. It can be seen that the total travel costs resulting from the stochastic MPC is less sensitive to the penetration rates than that of the multi-scale MPC. This is expected, as the proposed stochastic MPC controller aims to maintain the system at a desired state given that uncertain noises exist in the system. Regardless of the different levels of the noise, the control actions eventually brings the system towards the desired state. Furthermore, Notice that there are still noticeable gaps between the two curves even at a penetration rate of 100\%. One possible explanation is that the system has uncertainty in demand and MFD. Such uncertainty exists even for 100\%  penetration rate.

Next, let us look into more details at both the network level and the local level. At the network level (Fig.\ref{fig:cv}b), it is shown that the network travel time is not sensitive to the penetration rates for both the multi-scale and stochastic MPC. This is expected, as both controllers can stabilize the system. Note that the total network travel time for the scenarios with low penetration rates is even slightly smaller than that of the scenarios with high penetration rates. This is because there are less vehicles in the network for the scenarios with low penetration rates, as fewer vehicles are let in at the perimeter. At the local level (Fig.\ref{fig:cv}c), the total intersection delay is strongly influenced by the penetration rate for both the multi-scale and stochastic MPCs. By increasing the penetration rate from 5\% to 100\%, the total delay at the perimeter intersections resulting from the multi-scale MPC and the stochastic MPC is reduced by 41.1\% and 27.3\%, respectively. The reduction in delay clearly indicates that for both controllers, the penetration rate is more important for the intersection control than the network control. We can also see that the stochastic controller is less sensitive to the penetration rates than the multi-scale controller. This indicates that the stochastic controller exhibits more robustness to uncertainties in the system.

Even if the stochastic MPC tends to outperform the multi-scale MPC, the multi-scale MPC has the strong advantage of low computational cost. Hence, under higher penetration rates and low system stochasticity, it might be advisable to apply the multi-scale control in reality. 

\section{Conclusions and discussion}
This paper proposes a multi-scale MPC approach to integrate the network level perimeter control and the local level signal control in a traffic system with one center region and one periphery. 
The model calculates the optimal green ratios that minimize the total travel costs for both levels in a moving time horizon. 
Connected vehicles provide not only information on the current traffic accumulations in both network (traffic accumulation) and local level (queue length), but also provide  information on possible future arrivals. 
The predictive nature of the connected vehicles enables a more efficient application of the MPC approach.  
This MPC approach is further extended into a two-stage stochastic control scheme to cope with large system noises due to the stochasticity of the system and the lack of information. The stochastic controller optimizes the expected travel cost resulting from the current control action.  

Case studies show that the multi-scale MPC approach successfully stabilizes the traffic accumulation in the center network with the least impact on the perimeter intersections. 
The total travel cost is minimized for the vehicles both in the network and on the perimeter. 
It is also robust to moderate noises. 
In scenarios with large noises, the stochastic controller is shown to successfully reduce the oscillations in the system. 
We looked at the applicability of the proposed two controllers under different penetration rates of connected vehicles. 
For scenarios with very low penetration rate and/or large system noises, the stochastic controller is highly recommended. The sample size in the stochastic optimization scheme should be no less than 20 and compatible with the computational power of the system. For scenarios with high penetration rate and low system noises, multi-scale control already gives satisfactory results. 

In the proposed framework, we consider the inflow to the network entering from a fixed intersection and neglect the possible reroute to its adjacent intersections when queue length is large. Incorporating this effect will be a future work, which requires significant amount of research efforts on the modelling of the system dynamics. This work does not assume the outside network to be homogeneous with a single MFD. The arrival rate to the perimeter intersections is assumed to be predicted with real-time connected vehicles information. The impact of the control policy on the outside network and integration of the network and periphery for multi-region perimeter control are considered as future work. The model can also be adapted to address the spillbacks by incorporating either hard or soft spill back constraints.
\bibliography{reference}

\end{document}